\documentclass[12pt]{amsart}

\usepackage{times}
\usepackage{amssymb}
\usepackage{amsfonts}
\usepackage{amsthm}
\usepackage{bm}
\usepackage{graphicx}

\newcommand{\excise}[1]{}

\newtheorem{theorem}{Theorem}[section]
\newtheorem{lemma}[theorem]{Lemma}

\newtheorem{corollary}[theorem]{Corollary}
\newtheorem{prop}[theorem]{Proposition}

\newtheorem{questions}[theorem]{Questions}
\newtheorem{answers}[theorem]{Answers}

\theoremstyle{definition}
\newtheorem{example}[theorem]{Example}
\newtheorem{remark}[theorem]{Remark}
\newtheorem{definition}[theorem]{Definition}
\newtheorem{convention}[theorem]{Convention}
\newtheorem{observation}[theorem]{Observation}
\newtheorem{notation}[theorem]{Notation}

\numberwithin{equation}{section}

\newenvironment{numbered}%
        {\begin{list}
                {\noindent\makebox[0mm][r]{\arabic{enumi}.}}
                {\leftmargin=5.5ex \usecounter{enumi}}
        }
        {\end{list}}

        {\begin{list}
                {\noindent\makebox[0mm][r]{(\roman{enumi})}}
                {\leftmargin=5.5ex \usecounter{enumi}}
        }
        {\end{list}}

\newcounter{separated}

\newcommand{\ring}[1]{\ensuremath{\mathbb{#1}}}

\newcommand\0{\mathbf{0}}
\newcommand\<{\langle}
\renewcommand\>{\rangle}
\newcommand\CC{\ring{C}}
\newcommand\FF{{\mathcal F}}
\newcommand\HH{{\mathcal H}}
\newcommand\KK{{\mathcal K}}

\newcommand\NN{\ring{N}}
\newcommand\OO{{\mathcal O}}
\newcommand\QQ{\ring{Q}}
\newcommand\RR{\ring{R}}

\newcommand\VV{{\mathcal V}}
\newcommand\ZZ{\ring{Z}}

\newcommand\mm{{\mathfrak m}}
\newcommand\pp{{\mathfrak p}}

\newcommand\cB{{\mathcal B}}
\newcommand\cC{{\mathcal C}}
\newcommand\cI{{\mathcal I}}
\newcommand\cS{{\mathcal S}}
\newcommand\cZ{{\mathcal Z}}

\newcommand\oJ{{\hspace{.45ex}\overline{\hspace{-.45ex}J}}}

\newcommand\del{\partial}
\newcommand\ttt{\partial}
\newcommand\vea{\varepsilon_{\hspace{-.1ex}A}}

\newcommand\Horn{{\rm Horn}}
\newcommand\qdeg{{\rm qdeg}}
\newcommand\tdeg{{\rm tdeg}}

\renewcommand\th{{\rm th}}
\renewcommand\hom{\operatorname{Hom}}
\newcommand\ext{\operatorname{Ext}}
\newcommand\sat{{\rm sat}}
\newcommand\vol{{\rm vol}}
\newcommand\jump{\text{\rm jump}}
\newcommand\rank{{\rm rank}}
\newcommand\toral{\text{\rm toral}}
\newcommand\andean{\text{\rm Andean}}
\newcommand\primary{\text{\rm primary}}

\newcommand\onto{\twoheadrightarrow}
\newcommand\spot{{\hbox{\raisebox{1.7pt}{\large\bf .}}\hspace{-.5pt}}}
\newcommand\blank{\,\cdot\,}
\newcommand\minus{\smallsetminus}
\newcommand\nothing{\varnothing}

\def\ol#1{{\overline {#1}}}

\def\dlim{{\lim\limits_{\raisebox{.2ex}{$\scriptstyle\longrightarrow$}}}{_{\,}}}

\newcommand{\aoverb}[2]{{\genfrac{}{}{0pt}{1}{#1}{#2}}}
\def\twoline#1#2{\aoverb{\scriptstyle {#1}}{\scriptstyle {#2}}}

\begin{document}

\title{Binomial $D$-modules}

\author{Alicia Dickenstein}
\address[AD]{Departamento de Matem\'atica\\
FCEN, Universidad de Buenos Aires \\(1428) Buenos Aires, Argentina.}
\email{alidick@dm.uba.ar}
\thanks{AD was 
partially supported by UBACYT X042, CONICET PIP 5617 and ANPCyT PICT
20569, Argentina}

\author{Laura Felicia Matusevich}
\address[LFM]{Department of Mathematics \\
University of Pennsylvania \\ Philadelphia, PA 19104.}
\email{laura@math.tamu.edu}
\curraddr{Department of Mathematics, Texas A\&M University, College
Station, TX 77843.}
\thanks{LFM was partially supported by an NSF Postdoctoral Research
Fellowship}

\author{Ezra Miller}
\address[EM]{Department of Mathematics \\ University of Minnesota \\
Minneapolis, MN 55455.}
\email{ezra@math.umn.edu}
\thanks{EM was partially supported by NSF grants DMS-0304789 and DMS-0449102}

\dedicatory{The authors dedicate this work to the memory of Karin
Gatermann, friend and colleague\vspace{-4pt}}

\date{26 March 2008}

\begin{abstract}
We study quotients of the Weyl algebra by left ideals whose generators
consist of an arbitrary $\ZZ^d$-graded binomial ideal $I$ in
$\CC[\del_1,\ldots,\del_n]$ along with Euler operators defined by the
grading and a parameter $\beta \in \CC^d$.  We determine the
parameters~$\beta$ for which these $D$-modules (i)~are holonomic
(equivalently, regular holonomic, when $I$ is standard-graded);
(ii)~decompose as direct sums indexed by the primary components
of~$I$; and (iii)~have holonomic rank greater than the rank for
generic~$\beta$.  In each of these three cases, the parameters in
question are precisely those outside of a certain explicitly described
affine subspace arrangement in~$\CC^d$.  In the special case of Horn
hypergeometric $D$-modules, when $I$ is a lattice basis ideal, we
furthermore compute the generic holonomic rank combinatorially and
write down a basis of solutions in terms of associated
$A$-hypergeometric functions.  This study relies fundamentally on the
explicit lattice point description of the primary components of an
arbitrary binomial ideal in characteristic zero, which we derive in
our companion article \cite{primDecomp}.
\end{abstract}
\maketitle

\mbox{}
\vspace{-2.7ex}
\parskip=0ex
\parindent2em
\setcounter{tocdepth}{2}
\tableofcontents
\parskip=1ex
\parindent0pt

\section{Introduction}

\subsection{Hypergeometric series}\label{s:hypergeometric}

A univariate power series is \emph{hypergeometric}\/ if the successive
ratios of its coefficients are given by a fixed rational
function. These functions, and the elegant differential equations they
satisfy, have proven ubiquitous in mathematics. As a small example of
this phenomenon, consider the Hermite polynomials. These
hypergeometric functions naturally occur, for instance, in physics
(energy levels of the harmonic oscillator) \cite{quantum-mechanics}, 
numerical analysis (Gaussian quadrature) \cite{numerical-analysis-book},
combinatorics (matching polynomials of complete graphs) \cite{godsil},
and probability (iterated It\^{o} integrals of standard Wiener
processes) \cite{ito}.

Perhaps the most natural definition of hypergeometric power series in
several variables is the following, whose bivariate specialization was
studied by Jakob Horn as early as 1889 \cite{horn89}. More references
include \cite{horn31}, the first of six articles, all in Mathematische
Annalen between 1931 and 1940, and all containing ``Hypergeometrische
Funktionen zweier  Ver\"anderlichen'' (hypergeometric functions in two
variables) in their titles.  

\begin{definition}\label{def:hyp-power-series}
A formal series $F(z) = \sum_{\alpha \in \NN^m} a_\alpha
z_1^{\alpha_1} \cdots z_m^{\alpha_m}$ in $m$ variables with complex
coefficients is \emph{hypergeometric in the sense of Horn}\/ if there
exist rational functions $r_1, r_2,\ldots,r_m$ in $m$ variables such
that
\begin{equation}\label{eqn:def-hyp-series}
  \frac{a_{\alpha+e_k}}{a_\alpha}=r_k(\alpha) \quad \mbox{for all\;}
  \alpha \in \NN^m \; \mbox{and\;} k=1,\dots,m.
\end{equation}
Here we denote by $e_1,\ldots,e_m$ the standard basis vectors
of~$\NN^m$.
\end{definition}

Write the rational functions of the previous definition as
\[
  r_k(\alpha)=p_k(\alpha)/q_k(\alpha+e_k) \quad k=1,\dots m,
\]
where $p_k$ and $q_k$ are relatively prime polynomials.

Let $\del_{z_i}$ denote the partial derivative operator
$\frac{\del}{\del z_i}$.  Since for all monomial functions
$z^{\alpha}$ and polynomials $g$ we have
$g(z_1\del_{z_1},\ldots,z_m\del_{z_m}) z^\alpha =
g(\alpha_1,\ldots,\alpha_m) z^\alpha$, the series $F$ satisfies the
following \emph{Horn hypergeometric system of differential equations}:
\begin{equation}\label{eqn:horn1}
  q_k(z_1\del_{z_1},\dots, z_m\del_{z_m}) F(z) = z_k
  p_k(z_1\del_{z_1},\dots,z_m\del_{z_m}) F(z) \quad
  k=1,\dots,m,
\end{equation}
provided that, for $k=1,\dots,m$, the condition $q_k(\alpha)=0$ is
satisfied whenever $\alpha_k=0$.  See also
Remark~\ref{rem:solution-form}.

Of particular interest are the series where the numerators and
denominators of the rational functions $r_k$ factor into products of
linear factors.  (Contrast with the notion of ``proper hypergeometric
term'' in \cite{A=B}.)  Notice that by the fundamental theorem of
algebra, this is not restrictive when the number of variables is
$m=1$.

\subsection{Binomial ideals and binomial $D$-modules}\label{s:bin}

The central objects of study in this article are the \emph{binomial
$D$-modules}, to be introduced in Definition~\ref{d:Eulers}, which
reformulate and generalize the classical Horn hypergeometric systems,
as we shall see in Section~\ref{s:bihorn}.  Our definition is based on
the point of view developed by Gelfand, Graev, Kapranov, and
Zelevinsky \cite{GGZ, GKZ}, and contains their hypergeometric systems
as special cases; see Section~\ref{s:GKZ}.

To construct a binomial $D$-module, the starting point is an integer
matrix~$A$, about which we wish to be consistent throughout.

\begin{convention}\label{conv:A}
$A = (a_{ij}) \in \ZZ^{d \times n}$ denotes an integer $d \times n$
matrix of rank~$d$ whose columns $a_1,\ldots,a_n$ all lie in a single
open linear half-space of~$\RR^d$; equivalently, the cone generated by
the columns of $A$ is pointed (contains no lines), and all of the
$a_i$ are nonzero.  We also assume that $\ZZ A = \ZZ^d$; that is, the
columns of $A$ span $\ZZ^d$ as a lattice.
\end{convention}

The reformulation of Horn systems in Section~\ref{s:bihorn} proceeds
by a change of variables, so we will use $x = x_1,\ldots,x_n$ and
$\ttt = \del_1,\ldots,\del_n$ (where $\del_i = \del_{x_i} = \del/\del
x_i$), instead of $z_1,\ldots,z_m$ and $\del_{z_1},\ldots,\del_{z_m}$,
whenever we work in the binomial setting.  The matrix $A$ induces a
$\ZZ^d$-grading of the polynomial ring $\CC[\del_1,\dots,\del_n] =
\CC[\ttt]$, which we call the \emph{$A$-grading}, by setting
$\deg(\del_i) = -a_i$.  An ideal of $\CC[\ttt]$ is \emph{$A$-graded}\/
if it is generated by elements that are homogeneous for the
$A$-grading.  For example, a \emph{binomial ideal} is generated by
\emph{binomials} $\ttt^u - \lambda\ttt^v$, where $u,v \in \ZZ^n$ are
column vectors and $\lambda \in \CC$; such an ideal is $A$-graded
precisely when it is generated by binomials $\ttt^u - \lambda\ttt^v$
each of which satisfies either $Au = Av$ or $\lambda = 0$ (in
particular, monomials are allowed as generators of binomial ideals).
The hypotheses on~$A$ mean that the $A$-grading is a \emph{positive
$\ZZ^d$-grading} \cite[Chapter~8]{cca}.

The Weyl algebra $D = D_n$ of linear partial differential operators,
written with the variables $x$ and~$\ttt$, is also naturally
$A$-graded by additionally setting $\deg(x_i) = a_i$.  Consequently,
the \emph{Euler operators} in our next definition are $A$-homogeneous
of degree~$0$.

\begin{definition}\label{d:Eulers}
For each $i \in \{1,\ldots,d\}$, the $i^\th$ \emph{Euler operator}\/
is
\[
  E_i = a_{i1} x_1\del_1 + \cdots + a_{in} x_n\del_n.
\]
Given a vector $\beta \in \CC^d$, we write $E-\beta$ for the sequence
\mbox{$E_1 - \beta_1, \ldots, E_d - \beta_d$}.  (The dependence of the
Euler operators $E_i$ on the matrix~$A$ is suppressed from the
notation.)

For an $A$-graded binomial ideal $I \subseteq \CC[\ttt]$, we denote by
$H_A(I,\beta)$ the left ideal $I + \<E-\beta\>$ in the Weyl
algebra~$D$.  The \emph{binomial $D$-module} associated to~$I$ is
$D/H_A(I,\beta)$.
\end{definition}

We will explain in Section~\ref{s:bihorn} how Horn systems correspond
to the binomial $D$-modules arising from a very special class of
binomial ideals called \emph{lattice basis ideals}.

Our goal for the rest of this Introduction (and indeed, the rest of
the paper) is to demonstrate not merely that the definition of
binomial $D$-modules can be made in this generality---and that it
leads to meaningful theorems---but that it \emph{must} be made, even
if one is interested only in classical questions concerning
Horn hypergeometric systems, which arise from lattice basis ideals.
Furthermore, once the definition has been made, most of what we wish
to prove about Horn hypergeometric systems generalizes to all binomial
$D$-modules.

\subsection{Toric ideals and $A$-hypergeometric systems}\label{s:GKZ}

The fundamental examples of binomial $D$-modules, and the ones which
our definition most directly generalizes, are the
\emph{$A$-hypergeometric systems} (or \emph{GKZ hypergeometric
systems}) of Gelfand, Graev, Kapranov, and Zele\-vin\-sky
\cite{GGZ,GKZ}.  Given $A$ as in Convention~\ref{conv:A}, these are
the left $D$-ideals $H_A(I_A,\beta)$, also deno\-ted by $H_A(\beta)$,
where
\begin{equation}
\label{eq:IA} 
  I_A = \< \ttt^u -\ttt^v : Au=Av\> \subseteq \CC[\del_1,\dots,\del_n]
\end{equation}
is the \emph{toric ideal} for the matrix~$A$.  The systems
$H_A(\beta)$ have many applications; for example, they arise naturally
in the moduli theory of Calabi-Yau complete intersections in 
toric varieties, and (therefore)
they play an important role in applications of mirror symmetry in
mathematical physics \cite{bvs, horja, hosono, hly}.

The ideal $I_A$ is a prime $A$-graded binomial ideal, and the quotient
ring $\CC[\ttt]/I_A$ is the semigroup ring for the affine
semigroup~$\NN A$ generated by the columns of~$A$.  There is a rich
theory of toric ideals, toric varieties, and affine semigroup rings,
whose core philosophy is to exploit the connection between the algebra
of the semigroup ring $\CC[\ttt]/I_A = \CC[\NN A]$ and the
combinatorics of the semigroup~$\NN A$.  In this way,
algebro-geometric results on toric varieties can be obtained by
combinatorial means, and purely combinatorial facts about polyhedral
geometry can be proved using algebraic techniques.  We direct the
reader to the texts \cite{fulton-toric,gkz-book,cca} for more
information.
 
Much is known about $A$-hypergeometric $D$-modules.  They are
holonomic for all parameters \cite{GKZ,Adolphson}, and they are
regular holonomic exactly when $I_A$ is $\ZZ$-graded in the usual
sense \cite{equivariant, uli06}.  In this case, (Gamma-)series
expansions for the solutions of $H_A(\beta)$ centered at the origin
and convergent in certain domains can be explicitly computed
\cite{GKZ,SST}.  The generic (minimal) holonomic rank is known to be
$\vol(A)$, the normalized volume of the convex hull of the columns of
$A$ and the origin \cite{GKZ,Adolphson}, and holonomic rank is
independent of the parameter~$\beta$ if and only if the semigroup ring
$\CC[\NN A]$ is Cohen-Macaulay \cite{GKZ,Adolphson,MMW}.  We will
extend all of these results, suitably modified, to the general setting
of binomial $D$-modules.  The important caveat is that a general
binomial $D$-module can exhibit behavior that is forbidden to GKZ
systems (see Example~\ref{ex:non-holonomic!}, for instance), so it is
impossible for the extension to be entirely straightforward.

\subsection{Binomial Horn systems}\label{s:bihorn}

Classical Horn systems, which we are about to define precisely, were
first studied by Appell \cite{appell}, Mellin \cite{mellin}, and Horn
\cite{horn89}.  They directly generalize the univariate hypergeometric
equations for the functions $_{p}F_{q}$; see
\cite{multiple-gauss-book,slater} and the references therein.  As we
mentioned earlier, our motivation to consider binomial $D$-modules is
that they contain as special cases these classical Horn systems.  The
definition of these systems involves a matrix $B$ about which, like
the matrix $A$ from Convention~\ref{conv:A}, we wish to be consistent
throughout.

\begin{convention}\label{conv:B}
Let $B = (b_{jk})\in \ZZ^{n\times m}$ be an integer matrix of full
rank~$m \leq n$.  Assume that every nonzero element of the column-span
of~$B$ over the integers~$\ZZ$ is \emph{mixed}, meaning that it has at
least one positive and one negative entry; in particular, the columns
of $B$ are mixed.  We write $b_1,\ldots,b_n$ for the rows of~$B$.
Having chosen~$B$, we set $d = n - m$ and pick a matrix $A \in \ZZ^{d
\times n}$ whose columns span $\ZZ^d$ as a lattice, such that $AB =
0$. 

If $d\neq 0$, the mixedness hypothesis on $B$ is
equivalent to the pointedness assumption for~$A$ that appears in
Convention~\ref{conv:A}. We do allow~\mbox{$d=0$}, in which case $A$
is the empty matrix.
\end{convention}

\begin{definition}\label{def:old-horn}
For a matrix $B \in \ZZ^{n\times m}$ as in Convention~\ref{conv:B} and
a vector \mbox{$c=(c_1,\dots,c_n)$} in $\CC^n$, the \emph{classical
Horn system with parameter $c$}\/ is the left ideal $\Horn(B,c)$ in
the Weyl algebra~$D_m$ generated by the $m$ differential operators
\[
  q_k(\theta_z) - z_k p_k(\theta_z), \quad k= 1, \dots, m,
\]
where $\theta_z = (\theta_{z_1}, \ldots, \theta_{z_m})$, $\theta_{z_k}
= z_k \del_{z_k}$ $(1\leq k \leq m)$, and
\[
  q_k(\theta_z) =
  \prod_{b_{jk}>0}\prod_{\ell=0}^{b_{jk}-1}(b_j\cdot\theta_z+c_j-\ell)
  \quad \text{and} \quad p_k(\theta_z) =
  \prod_{b_{jk}<0}\prod_{\ell=0}^{|b_{jk}|-1}(b_j\cdot\theta_z+c_j-\ell).
\]
\end{definition}

\begin{remark}\label{rem:solution-form}
When the parameter $c$ is generic, one can find a
local basis of solutions of $\Horn(B,c)$ that consists of Puiseux
series of the form $z^vF(z)$, for certain complex vectors $v$ and
power series $F$ that are hypergeometric in the sense of
Definition~\ref{def:hyp-power-series}.  The rational functions giving
the recursions for the coefficients of these series~$F$ are related to
the defining equations for $\Horn(B,c)$.  We can see this more clearly
in an example.
\end{remark}

\begin{example}\label{ex:gauss-hyp-equation}
For the column matrix
$
B = \left[\begin{array}{cccc} 1 & -1 & -1 & 1 \end{array} \right]^t,
$
the corresponding Horn system with parameter $(c_1,c_2,c_3,c_4)$
consists of one operator in the single variable $z$, namely
\[
(\theta_z + c_1)(\theta_z + c_4) - z (-\theta_z + c_2)(-\theta_z + c_3)
=
(\theta_z + c_1)(\theta_z + c_4) - z (\theta_z - c_2)(\theta_z - c_3)
.
\]
We can follow the usual convention of normalizing $c_4$ to $1$ and
renaming the parameters to obtain the operator
\begin{equation}\label{eqn:gauss}
  (\theta_z + c)(\theta_z + 1) - z(\theta_z+a)(\theta_z+b), \quad
  (\mbox{here } a,b,c \in \CC).
\end{equation}
This is the Gauss hypergeometric equation multiplied on the left by
the variable $z$ (this does not alter the space of local holomorphic
solutions), and written in operator notation.

If $c$ is not an integer, we can write down a local basis of solutions
for~(\ref{eqn:gauss}) converging in a disk centered at the origin.
This basis consists of the functions $F(z)$ and $z^{1-c}G(z)$, where
\[
  F(z) = 1 + \frac{ab}{1!\,c} z + \frac{a(a+1)b(b+1)}{2!\,c(c+1)} z^2+
  \cdots,
\]
\[
\newcommand\+{\hspace{-.3ex}+\hspace{-.3ex}}
\newcommand\?{\hspace{-.3ex}-\hspace{-.3ex}}
  G(z) = 1 + \frac{(a\+1\?c)(b\+1\?c)}{1!\,(2\?c)} z +
  \frac{(a\+1\?c)(a\+1\?c\+1)(b\+1\?c)(b\+1\?c\+1)}
  {2!\,(2\?c)(2\?c\+1)} z^2 + \cdots\!.
\]
The rational functions giving the recursions for the
coefficients of the hypergeometric series $F$ and $G$ are
\[
  r(\alpha) = \frac{(\alpha+a)(\alpha+b)}{(\alpha+c)(\alpha+1)}
  \quad\text{and}\quad
  s(\alpha) = \frac{(\alpha+(a+1-c))(\alpha+(b+1-c))}{(\alpha+(2-c))(\alpha+1)},
\]
respectively; the numerators and denominators of these rational
functions bear a marked (and non-coincidental) resemblance to the
polynomials in $\theta_z$ that appear in the two terms
of~(\ref{eqn:gauss}).
\end{example}

Using ideas of Gelfand, Kapranov, and Zelevinsky, the classical Horn
systems can be reinterpreted as the following binomial $D$-modules,
with $\beta = Ac$.

\begin{definition}\label{def:new-horn}
Fix integer matrices $B$ and $A$ as in Convention~\ref{conv:B}, and
let $I(B)$ be the \emph{lattice basis ideal} corresponding to this
matrix, that is, the ideal in $\CC[\ttt]$ generated by the
binomials
\[
  \prod_{b_{jk}>0} \del_{x_j}^{b_{jk}} - \prod_{b_{jk}<0}
  \del_{x_j}^{-b_{jk}} \quad\text{for}\quad 1 \leq k \leq m.
\]
The \emph{binomial Horn system with parameter~$\beta$}\/ is the left
ideal $H(B,\beta) = H_A(I(B),\beta)$ in the Weyl algebra $D = D_n$.
\end{definition}

The classical-to-binomial transformation proceeds via the surjection
\begin{equation}\label{eqn:change-of-vars}
  \begin{array}{rcl}
         (\CC^*)^n &\rightarrow& (\CC^*)^m \\
  (x_1,\ldots,x_n) &  \mapsto  & x^B =
  (\prod_{j=1}^{n}x_j^{b_{j1}},\ldots,\prod_{j=1}^nx_j^{b_{jm}}),
  \end{array}
\end{equation}
where $\CC^* = \CC \minus \{0\}$ is the group of nonzero complex
numbers.  A solution $f(z_1,\dots,z_m)$ of the classical Horn system
$\Horn(B,c)$ gives rise to a solution~$x^cf(x^B)$ of the binomial Horn
system~$H(B,Ac)$.  When the columns of $B$ are a basis of the integer kernel
of $A$, this map  defines a vector space isomorphism
between the (local) solution spaces. This was proved in
\cite[Section~5]{dms} for $n > m$ in the homogeneous case, where the
column sums of~$B$ are zero, but the proofs (which are elementary
calculations taking only a page) go through verbatim for $n\geq m$ in
the inhomogeneous case.

The transformation $f(z) \mapsto x^cf(x^B)$ takes classical series
solutions supported on $\NN^m$ to Puiseux series solutions supported
on the translate $c + \ker(A) \subseteq \CC^n$ of the kernel of~$A$
in~$\ZZ^n$.  (Note that $\ker(A)$ contains the lattice $\ZZ B$ spanned
by the columns of~$B$ as a finite index subgroup.)  More precisely,
the differential equations $E-\beta$, which geometrically impose
torus-equivariance infinitesimally under the action of (the Lie
algebra of) \mbox{$\ker((\CC^*)^n \to (\CC^*)^m)$}, result in series
supported on $c + \ker(A)$, while the binomials in the lattice basis
ideal $I(B) \subseteq H(B,Ac)$ impose hypergeometric constraints on
the coefficients.

Although the isomorphism $f(z) \mapsto x^c f(x^B)$ is only at the
level of local holomorphic solutions, not $D$-modules, it preserves
many of the pertinent features, including the dimensions of the spaces
of local holomorphic solutions and the structure of their series
expansions.  Therefore, although the classical Horn systems are our
motivation, we take the binomial formulation as our starting point: no
result in this article depends logically on the classical-to-binomial
equivalence.

\subsection{Holomorphic solutions to Horn systems}\label{s:sols}

The binomial rephrasing of Horn systems led to formulas in \cite{GGR}
for Gamma-series solutions via $A$-hypergeometric theory.  However,
Gamma-series need not span the space of local holomorphic solutions
of~$H(B,\beta)$ at a point of~$\CC^n$ that is nonsingular for
$H(B,\beta)$, even in the simplest cases.  The reason is that
Gamma-series are \emph{fully supported}: there is a cone of
dimension~$m$ (the maximum possible) whose lattice points correspond
to monomials with nonzero coefficients.  Generally speaking, Horn
systems in dimension $m \geq 2$ tend to have many series solutions
\mbox{without full support}.

\begin{example}\label{ex:erdelyi}
In the course of studying one of Appell's systems of two
hypergeometric equations in $m = 2$ variables, Arthur Erd\'elyi
\cite{erdelyi} mentions a modified form of the following example.
Given any $\beta \in \CC^2$ and the two matrices
\[
A =   \left[ \begin{array}{cccc}
       3 & 2 & 1 & 0 \\
       0 & 1 & 2 & 3 
\end{array} \right], \quad
B =  \left[ \begin{array}{rr}
       1  & 0  \\
       -2 & 1  \\
       1  & -2 \\
       0  &  1 
\end{array} \right]
\]
satisfying Convention~\ref{conv:B}, the Puiseux monomial
$x_1^{\beta_1/3}x_4^{\beta_2/3}$ is a solution of~$H(B,\beta)$.
\end{example}

A key feature of the above example is that the solutions without full
support persist for arbitrary choices of the parameter vector~$\beta$.
The fact that this phenomenon occurs in much more generality---for
arbitrary dimension $m \geq 2$, in particular---was realized only
recently \cite{dms}.  And it is not the sole peculiarity that arises
in dimension $m \geq 2$: in view of the transformation to binomial
Horn systems in Section~\ref{s:bihorn}, the following demonstrates
that classical Horn systems can exhibit poor behavior for badly chosen
parameters.

\begin{example}\label{ex:non-holonomic!}
Consider 
\[
  A = \left[ \begin{array}{rrrr}
                1 & 1 & 0 & 0\\
                0 &  1 & 1 & 1
             \end{array}
      \right]
  \quad\text{and}\quad
  B = \left[ \begin{array}{rr}
                1 & 1 \\
               -1 & -1\\
                1 & 0 \\ 
                0 & 1 
             \end{array}
      \right],
\]
so that
\[
  H(B,\beta) = \< \del_1\del_3-\del_2,\del_1\del_4-\del_2\> + 
  \<x_1\del_1-x_2\del_2-\beta_1, x_2\del_2+x_3\del_3+x_4\del_4-\beta_2\>.
\]
If $\beta_1=0$, then any (local holomorphic) bivariate function
$f(x_3,x_4)$ annihilated by the operator $x_3\del_3+x_4\del_4-\beta_2$
is a solution of $H(B,\beta)$.  The space of such functions is
infinite-dimensional; in fact, it has uncountable dimension, as it
contains all monomials $x_3^{w_3} x_4^{w_4}$ with $w_3,w_4 \in \CC$
and $w_3+w_4 = \beta_2$.
\end{example}

Erd\'elyi's goal for his study of the Appell system was to give bases
of solutions that converged in different regions of~$\CC^2$,
eventually covering the whole space, just as Kummer had done for the
Gauss hypergeometric equation more than a century before
\cite{kummer}.  There has been extensive work since then (see
\cite{multiple-gauss-book} and its references) on convergence of more
general hypergeometric functions in two and three variables.  But
already for the classical case of Horn systems, where the phenomena in
Examples~\ref{ex:erdelyi} and~\ref{ex:non-holonomic!} are commonplace,
Erd\'elyi's work raises a number of fundamental questions that 
have remained largely open (partial answers in dimension $m = 2$ being
known \cite{dms}; see Remark~\ref{rk:dms}).  The purpose of this
article is to answer the following completely and precisely.

\begin{questions}\label{q}
Fix~$B$ as in Convention~\ref{conv:B} and consider the Horn systems
determined~by~$B$.
\begin{numbered}
\item\label{q:finite-rank}
For which parameters does the space of local holomorphic solutions
around a nonsingular point have finite dimension as a complex vector
space?

\item\label{q:rank}
What is a combinatorial formula for the minimum such dimension, over
all possible choices of parameters?

\item\label{q:generic}
Which parameters are generic, in the sense that the minimum dimension
is attained?

\item\label{q:support}
How do (the supports of)  series solutions centered at the origin
look, combinatorially?
\end{numbered}\setcounter{separated}{\value{enumi}}%
\end{questions}

These questions make sense simultaneously for classical Horn systems
and binomial Horn systems, since the answers are invariant under the
classical-to-binomial transformation.  That the questions also make
sense for binomial $D$-modules is our point of departure, for they can
be addressed in this generality using answers to the following.

\addtocounter{theorem}{-1}
\begin{questions}[continued]
Consider the binomial $D$-modules $H_A(I,\beta)$ for varying
\mbox{$\beta \!\in\! \CC^d$}.
\begin{numbered}\setcounter{enumi}{\value{separated}}
\item\label{q:hol}
When is $D/H_A(I,\beta)$ a holonomic $D$-module? 

\item\label{q:reg}
When is $D/H_A(I,\beta)$ a \emph{regular} holonomic $D$-module?
\end{numbered}
\end{questions}

The phenomena underlying all of the answers to Questions~\ref{q} can
be described in terms of lattice point geometry, as one might hope,
owing to the nature of hypergeometric recursions as relations between
coefficients on monomials.  The lattice point geometry is elementary,
in the sense that it only requires constructions involving cosets and
equivalence relations in lattices.  However, modern techniques are
required to make the descriptions quantitatively accurate and prove
them.  In particular, our progress applies two distinct and
substantial steps: precise advances in the combinatorial commutative
algebra of binomial ideals in semigroup rings \cite{primDecomp}, and
the functorial translation of those advances into $D$-module theory
here.

\subsection{Combinatorial answers to hypergeometric questions}\label{s:latgeom}

The supports of the various series solutions to~$H(B,\beta)$ centered
at the origin are controlled by how effectively the columns of~$B$
join the lattice points in the positive orthant~$\NN^n$.  In essence,
this is because the coefficients on a pair of Puiseux monomials are
related by the binomial equations in~$I(B)$ when their exponent
vectors in $c + \ker(A)$ differ by a column of~$B$.  This observation
prompts us to construct an undirected graph $\Gamma_B(\NN^n)$ on the
nodes~$\NN^n$ with an edge between pairs of points differing by a
column of~$B$.  Each connected component, or \emph{$B$-subgraph
of\/~$\NN^n$}, is contained in a single fiber of the projection $\NN^n
\to \ZZ^n/\ZZ B$.

Certain pairs consisting of a subset $J \subseteq \{1,\ldots,n\}$ and
a saturated sublattice $L \subseteq \ZZ^J$ contained in $\ker(A)$ are
\emph{associated}\/ to~$B$.  The columns of~$B$ together with the
vectors in~$L$ determine a graph $\Gamma_{B,L}(\ZZ^J \times \NN^\oJ)$
with nodes $\ZZ^J \times \NN^\oJ$, where $\oJ$ is the complementary
subset.  Each connected component of $\Gamma(B,L)$ is acted upon
by~$L$ and hence is a union of cosets of~$L$.  The key feature of an
associated lattice $L \subseteq \ZZ^J$ is that some of these
components consist of only finitely many cosets of~$L$; let us call
these components \emph{$L$-bounded}.

The $D$-module theoretic consequences of $L$-bounded components rely
on a crucial distinction; see Definition~\ref{d:toral},
Definition~\ref{d:andean}, and Remark~\ref{rk:andean} for more
precision and an etymology.

\begin{definition}\label{def:prelim-toralAndean}
An associated saturated sublattice $L \subseteq \ZZ^J \cap \ker(A)$ is
called \emph{toral}\/ if $L = \ZZ^J \cap \ker(A)$; otherwise, $L
\subsetneq \ZZ^J \cap \ker(A)$ is called \emph{Andean}.
\end{definition}

\begin{example}{}[Example~\ref{ex:non-holonomic!} continued]\label{ex:toralAndean1} 
With $A$ and $B$ as in Example~\ref{ex:non-holonomic!},
there are two associated lattices, one with $J=\{1,2,3,4\}$,
the other with $J = \{ 3,4\}$. The first one is toral,
while the second is Andean.
\end{example}

In what follows, $A_J$ denotes the submatrix of~$A$ whose columns are
indexed by~$J$.  We write $\ZZ A_J \subseteq \ZZ^d = \ZZ A$ for the
group generated by these columns, and $\CC A_J \subseteq \CC^d$ for
the vector subspace they generate.

\begin{observation}[{cf.\ \cite[Theorem~3.2]{primDecomp} and Lemma~\ref{l:arrangement}}]
\label{o:andean}
The images in~$\ZZ A$
of the $L$-bounded components for all of the Andean associated
sublattices $L \subseteq \ZZ^J$ comprise a finite union of cosets
of~$\ZZ A_J$.  The union over all~$J$ of the corresponding cosets of\/
$\CC A_J$ is an affine subspace arrangement in~$\CC^d$ called the
\emph{Andean arrangement} (Definition~\ref{d:Z} and~Lemma~\ref{l:Z}).
\end{observation}

\begin{example}{}[Example~\ref{ex:toralAndean1} continued]
The Andean arrangement in this case is 
\[
\CC \left[ \begin{smallmatrix}
0 & 0 \\ 1 & 1
\end{smallmatrix}
\right] =
\big\{ \left[ \begin{smallmatrix} 
0 \\ \beta_2 
\end{smallmatrix} \right]: \beta_2\in \CC\big\}.
\]
As we have already checked,
the Horn system in Example~\ref{ex:non-holonomic!}
fails to be holonomic for this set of parameters.
\end{example}

\begin{observation}[{cf.\ \cite[Theorem~4.12]{primDecomp} and its proof}]
\label{o:toral}
A component in $\ZZ^J \times \NN^\oJ$ determined by a toral associated
sublattice $L \subseteq \ZZ^J$ is $L$-bounded if and only if its image
in $\NN^\oJ$ is bounded.  If $\CC A_J = \CC^d$, the number of such
bounded images in $\NN^\oJ$ is finite; let $\mu(L,J)$ be the product
of this number with 
the index $|L/(\ZZ B \cap \ZZ^J)|$ of the sublattice~\mbox{$\ZZ B \cap
\ZZ^J$~in~$L$}.
\end{observation}

\begin{answers}\label{a}
The answers to Questions~\ref{q}, phrased in the language of binomial
Horn systems $H(B,\beta)$, are as follows.
\begin{numbered}
\item (Theorem~\ref{t:holonomic})
The dimension is finite exactly for $-\beta$ not in the Andean
arrangement.

\item (Theorem~\ref{t:Irank})
The generic (minimum) rank is $\sum \mu(L,J) \cdot \vol(A_J)$, the sum
being over all toral associated sublattices with $\CC A_J = \CC^d$,
where $\vol(A_J)$ is the volume of the convex hull of~$A_J$ and the
origin, normalized so a lattice simplex in~$\ZZ A_J$ has~volume~$1$.

\item (Definition~\ref{d:jump} and Theorem~\ref{t:Irank})
The minimum rank is attained precisely when $-\beta$ lies outside of
an affine subspace arrangement determined by certain local cohomology
modules, with the same flavor as (and containing) the Andean
arrangement.

\item (Theorem~\ref{t:Irank}, Theorem~\ref{thm:solsfromcomp}, and
Corollary~\ref{cor:supports})
\label{a:support}
When the Horn system is regular holonomic and~$\beta$ is general,
there are $\mu(L,J) \cdot \vol(A_J)$ linearly independent Puiseux
series solutions supported on (translates of) $L$-bounded components,
with coefficients determined by hypergeometric recursions.  Only $g \cdot
\vol(A)$ many Gamma-series solutions have full support, where $g =
|\ker(A)/\ZZ B|$ is the index of\/~$\ZZ B$ in its saturation.

\item (Theorem~\ref{t:holonomic})
Holonomicity is equivalent to the finite dimension in
Answer~\ref{a}.\ref{q:finite-rank}.

\item (Theorem~\ref{t:holonomic})
Holonomicity is equivalent to regular holonomicity when $I$ is
standard $\ZZ$-graded---i.e., the row-span of~$A$ contains the vector
$[1 \cdots 1]$.  Conversely, if there exists a parameter $\beta$ for
which $D/H_A(I,\beta)$ is regular holonomic, then $I$ is $\ZZ$-graded.
\end{numbered}
\end{answers}

In Answer~\ref{a}.\ref{a:support}, the solutions for toral sublattices
$L = \ker(A) \cap \ZZ^J$ in which $J$ is a proper subset of
$\{1,\ldots,n\}$ give rise to solutions that are bounded in the
$\NN^\oJ$ directions, and hence supported on sets of dimension
$\rank(L) = |J|-d < n-d = m$.  Answer~\ref{a}.\ref{q:reg} is, given
the other results in this paper, an (easy) consequence of the (hard)
holonomic regularity results of Hotta \cite{equivariant} and
Schulze--Walther \cite{uli06}.  Finally, let us note again that most
of the theorems quoted in Answers~\ref{a} are stated and proved in the
context of arbitrary binomial $D$-modules, not just Horn systems.

\begin{remark}\label{rk:dms}
We concentrate on the special case of Horn systems in
Section~\ref{s:horn}.  The systematic study of binomial Horn systems
was started in \cite{dms} under the hypothesis that $m$ (the number of
columns of $B$) is equal to~2. See also \cite{timur}.
Our results here are more general than
those found in \cite{dms} (as we treat all binomial $D$-modules, not
just those arising from lattice basis ideals of codimension~2), more
refined (we have completely explicit control over the parameters) and
stronger (for instance, our direct sum results hold at the level of
$D$-modules and not just local solution spaces).  On the other hand,
the generic holonomicity of classical Horn $D$-modules
(Definition~\ref{def:old-horn}) for $m>2$ remains unproven, the
bivariate case having been treated in \cite{dms}.
\end{remark}

\begin{example}{}[Example~\ref{ex:erdelyi}, continued]\label{ex:erdelyi'}
There are two associated sublattices $L \subseteq \ZZ^J$ here, both
toral, and both satisfying $\CC A_J = \CC^2$: the sublattice $\ker(A)
\subseteq \ZZ^4$, where $J = \{1,2,3,4\}$, and the sublattice $\0
\subseteq \ZZ^J$ for $J = \{1,4\}$.  Both of the multiplicities
$\mu\big(\ker(A),\{1,2,3,4\}\big)$ and $\mu\big(\0,\{1,4\}\big)$
equal~$1$, while $\vol(A) = 3$ and $\vol(A_{\{1,4\}}) = 1$, the latter
because $\left[\twoline 30\right]$ and $\left[\twoline 03\right]$ form
a basis for the lattice they generate.  Hence there are four solutions
in total, three of them with full support and one---namely the Puiseux
monomial in Example~\ref{ex:erdelyi}---with support of dimension zero.
See Example~\ref{ex:erdelyi''} for an (easy!) computation of these
associated lattices and their multiplicities.
\end{example}

\subsection{Binomial primary decomposition}\label{s:binom}

Our combinatorial study of binomial primary decomposition in
\cite{primDecomp} results in a natural language for quantifying which
sublattices are associated, which cosets appear in
Observation~\ref{o:andean}, and which bounded images appear in
Observation~\ref{o:toral}.  To be precise, a binomial prime ideal
$I_{\rho,J}$ in~$\CC[\del_1,\ldots,\del_n]$ is determined by a subset
$J \subseteq \{1,\ldots,n\}$ and a character $\rho: L \to \CC^*$ for
some sublattice $L \subseteq \ZZ^J$.  The sublattice $L \subseteq
\ZZ^J$ is associated to~$I(B)$, in the language of
Section~\ref{s:latgeom}, when $I_{\rho,J}$ is associated to~$I(B)$ in
the usual commutative algebra sense, and the multiplicity $\mu(L,J)$
of~$L$ in $I(B)$ from Observ\-ation~\ref{o:toral} is
$|L/(\ZZ B\cap\ZZ^J)|$ times the commutative algebra multiplicity
of~$I_{\rho,J}$~in~$I(B)$.  The factor of $|L/(\ZZ B\cap\ZZ^J)|$
counts the number of partial characters $\rho : L \to \CC^*$ for which
$I_{\rho,J}$ is associated to~$I(B)$; see Remark~\ref{rk:I(B)} for the
general reason~why.

\begin{example}{}[Example~\ref{ex:erdelyi'}, continued]\label{ex:erdelyi''}
The binomial Horn system is
\[
  H(B,\beta) = I(B) +  
       \<3 x_1\del_1 + 2 x_2\del_2 +  x_3\del_3              - \beta_1,
                         x_2\del_2 + 2x_3\del_3 + 3x_4\del_4 - \beta_2\>
  \subseteq D_4.
\]
The primary decomposition of the lattice basis ideal~$I(B)$ in
$\CC[\del_1,\del_2,\del_3,\del_4]$ is
\[
  I(B) = \<\del_1\del_3 - \del_2^2, \del_2\del_4 - \del_3^2\> =
  \<\del_1\del_3 - \del_2^2, \del_2\del_4 - \del_3^2, \del_1\del_4 -
  \del_2\del_3\> \cap \<\del_2, \del_3\>.
\]
The first of these components is the toric ideal $I_A = I_{\rho,J}$ of
the twisted cubic curve, where $\rho: \ker(A) = \ZZ B \to \CC^*$ is
the trivial character and $J = \{1,2,3,4\}$.  The ideal
$\<\del_2,\del_3\>$ is the binomial prime ideal $I_{\rho,J}$ for the
(automatically) trivial character $\rho: \0 \to \CC^*$ and the subset
$J = \{1,4\}$.  Both of these ideals have multiplicity~$1$ in~$I(B)$,
which is a radical ideal.  This explains the associated lattices and
multiplicities in Example~\ref{ex:erdelyi'}.
\end{example}

For a note on motivation, this project began with the conjectural
statement of \mbox{Theorem~\ref{thm:solsfromcomp}}
(Answer~\ref{a}.\ref{q:support}), which we concluded must hold because
of evidence derived from our knowledge of series solutions.  Its proof
reduced quickly to the statement of Example~\ref{ex:I(B)'}, which
directed all of the developments in the rest of the paper and in
\cite{primDecomp}.  Our consequent application of $B$-sub\-graphs and
their generalizations toward the primary decomposition of binomial
ideals serves as an advertisement for hypergeometric intuition as
inspiration for developments of independent interest in combinatorics
and commutative algebra.

\subsection{Euler-Koszul homology}\label{s:introEK}

Binomial primary decomposition is not only the natural language for
lattice point geometry, it is the reason why lattice point geometry
governs the $D$-module theoretic properties of binomial $D$-modules.
This we demonstrate by functorially translating the commutative
algebra of $A$-graded primary decomposition directly into the
$D$-module setting.  The functor we employ is Euler-Koszul homology
(see the opening of Section~\ref{s:EK} for background and references),
which allows us to pull apart the primary components of binomial
ideals, thereby isolating the contribution of each to the solutions of
the corresponding binomial $D$-module.  Here we see again the need to
work with general binomial $D$-modules: primary components of lattice
basis ideals, and intersections of various collections of them, are
more or less arbitrary $A$-homogeneous binomial ideals.

We stress at this point that the combinatorial geometric lattice-point
description of binomial primary decomposition is a crucial
prerequisite for the effective translation into the realm of
$D$-modules.  Indeed, semigroup gradings pervade the arguments
demonstrating the fundamentally holonomic behavior of Euler-Koszul
homology for toral modules (Theorem~\ref{t:rigid}) and its resolutely
non-holonomic behavior for Andean modules (Corollary~\ref{c:rigid}).
This is borne out in Lemma~\ref{l:toralAndean} and
Example~\ref{ex:andean}, which say that quotients by binomial primary
ideals are either toral or Andean as $\CC[\ttt]$-modules, thus
constituting the bridge from the commutative binomial theory in
\cite{primDecomp} to the binomial $D$-module theory in
Sections~\ref{s:EK} and~\ref{s:binomial}.  Taming the homological
(holonomic) and structural properties of binomial $D$-modules in
Theorems~\ref{t:holonomic}, \ref{t:sum}, and~\ref{t:Irank}---which,
together with Theorem~\ref{thm:solsfromcomp} on series bases, form our
core results---also rests squarely on having tight control over the
interactions of primary decomposition with various semigroup gradings
of the polynomial ring.  The underlying phenomenon is thus:

\textbf{Central principle.}
Just as toric ideals are the building blocks of binomial ideals,
$A$-hyper\-geometric systems are the building blocks of binomial
$D$-modules.

As a final indication of how structural results for binomial
$D$-modules have concrete combinatorial implications for Horn
hypergeometric systems, let us see how the primary decomposition in
Example~\ref{ex:erdelyi''} results in the combinatorial multiplicity
formula (Answer~\ref{a}.\ref{q:rank}) for the holonomic rank at
generic parameters~$\beta$.  The general result to which we appeal is
Theorem~\ref{t:sum}: for generic parameters~$\beta$, the binomial
$D$-module $D/H_A(I,\beta)$ decomposes as a direct sum over the toral
primary components of~$I$.

\begin{example}{}[Example~\ref{ex:erdelyi''}, continued]\label{ex:erdelyi'''}
The intersection in $\CC^4 = \mathrm{Spec}(\CC[\del_1,\dots,\del_4])$ of the two
irreducible varieties in the zero set of~$I(B)$ is the zero set of
\[
 \<\del_1\del_3 - \del_2^2, \del_2\del_4 - \del_3^2, \del_1\del_4 -
 \del_2\del_3\> + \<\del_2, \del_3\> = \<\del_1\del_4,\del_2,\del_3\>.
\]
The primary arrangement in Theorem~\ref{t:sum} is, in this case, the
line in~$\CC^2$ spanned by $\left[\twoline 30\right]$ union the line
in~$\CC^2$ spanned by $\left[\twoline 03\right]$.  When $\beta$ lies
off the union of these two lines, Theorem~\ref{t:sum} yields an
isomorphism of $D_4$-modules:
\[
  \frac{D_4}{H(B,\beta)} \cong 
  \frac{D_4}{\<\del_1\del_3 - \del_2^2, \del_2\del_4 - \del_3^2,
  \del_1\del_4 - \del_2\del_3\> + \<E-\beta\>} \oplus
  \frac{D_4}{\<\del_2,\del_3\> + \<E-\beta\>}.
\]
The summands on the right-hand side are GKZ hypergeometric systems (up
to extraneous vanishing variables in the $\<\del_2,\del_3\>$ case)
with holonomic ranks~$3$ and~$1$, respectively.
\end{example} 

We conclude this introduction with some general background on
$D$-modules.  A left $D$-ideal $\cI$ is \emph{holonomic}\/ if its
characteristic variety has dimension $n$.  Holonomicity has strong
homological implications, making the class of holonomic $D$-modules a
natural one to study. If $\cI$ is holonomic, its \emph{holonomic
rank}, i.e.\ the dimension of the space of solutions of the $D$-ideal
$\cI$ that are holomorphic in a sufficiently small neighborhood of a
point outside the singular locus, is finite (the converse of this
result is not true).  We refer to the texts \cite{borel, coutinho,SST}
for introductory overviews of the theory of $D$-modules; we point out
that the exposition in \cite{SST} is geared toward algorithms and
computations.  A treatment of $D$-modules with \emph{regular
singularities}\/ can be found in \cite{Bjork, Bjork2}.

\subsection*{Acknowledgments}

We wish to thank Uli Walther and Yves Laurent for inspiring
conversations.  This project benefited greatly from visits by its
authors to the University of Pennsylvania, Texas A\&M University, the
Institute for Mathematics and its Applications (IMA) in Minneapolis,
the University of Minnesota, and the Centre International de
Rencontres Math\'ematiques in Luminy (CIRM).  We thank these
institutions for their \mbox{gracious hospitality}.

\section{Euler-Koszul homology}\label{s:EK}

The Euler operators in Definition~\ref{d:Eulers} can be used to build
a Koszul-like complex whose zeroth homology is the $A$-hypergeometric
system in Definition~\ref{d:Eulers}.  In its most basic form, this
construction is due to Gelfand, Kapranov, and Zelevinsky \cite{GKZ},
and was developed by Adolphson \cite{Adolphson,Adolphson-Rend} and
Okuyama \cite{okuyama}, among others.  A~functorial generalization was
introduced in \cite{MMW}, where it was proved to be
homology-isomorphic to an ordinary Koszul complex detecting holonomic
rank changes for varying parameters~$\beta$.  Here we review the
definitions from \cite[Section~4]{MMW} (where more details can be
found), as well as connection to the quasidegrees defined in
\cite[Section~5]{MMW}.

Given a matrix~$A$ with columns $a_1,\dots,a_n$ as in
Convention~\ref{conv:A}, recall that the polynomial ring $\CC[\ttt]$
and the Weyl algebra $D=D_n$ are $A$-graded by $\deg(\del_j) = -a_j$
and $\deg(x_j)=a_j$.  Under this $A$-grading, operators
$E_1,\ldots,E_d$, and in fact all of the products $x_j\del_j\in D$,
are homogeneous of degree~$0$.

Given an $A$-graded left $D$-module~$W$, if~\mbox{$z \in W_\alpha$} is
homogeneous of degree~$\alpha$ then set $\deg_i(z) = \alpha_i$.  The
map $E_i - \beta_i: W \to W$ that sends each homogeneous element $z
\in W$ to
\begin{equation}\label{eq:E}
  (E_i - \beta_i) \circ z = (E_i-\beta_i-\deg_i(z)) z,
\end{equation}
and is extended $\CC$-linearly to all of~$W$, determines a $D$-linear
endomorphism of~$W$.

\begin{definition}\label{d:EK}
Fix $\beta\in\CC^d$ and an $A$-graded $\CC[\ttt]$-module~$V$.  The
\emph{Euler-Koszul complex}\/ $\KK_\spot(E-\beta;V)$ is the Koszul
complex of left $D$-modules defined by the sequence $E-\beta$ of
commuting endomorphisms on the left $D$-module $D\otimes_{\CC[\ttt]}
V$ concentrated in homological degrees $d$ to $0$.  The $i^\th$
\emph{Euler-Koszul homology}\/ of~$V$ is
\mbox{$\HH_i(E-\beta;V)=H_i(\KK_\spot(E-\beta;V))$}.
\end{definition}

\begin{example}\label{ex:A-hyp}
Fix $A$ and~$B$ as in Conventions~\ref{conv:A} and~\ref{conv:B}.
\begin{numbered}
\item
The \emph{binomial Horn $D$-module}\/ with parameter $\beta$ is
$\HH_0(E-\beta;\CC[\ttt]/I(B))$.
\item
The \emph{$A$-hypergeometric $D$-module}\/ with parameter $\beta$ is
$\HH_0(E-\beta;\CC[\ttt]/I_A)$; see~(\ref{eq:IA}).
\item
If $I \subseteq \CC[\ttt]$ is any $A$-graded binomial ideal, then
$\HH_0(E-\beta;\CC[\ttt]/I) = D/H_A(I,\beta)$.
\end{numbered}
\end{example}

Euler-Koszul homology behaves predictably with regard to $A$-graded
translation.

\begin{lemma}\label{l:V}
Let $V$ be an $A$-graded $\CC[\ttt]$-module and $\alpha \in \ZZ^d=\ZZ
A$.  If $V(\alpha)$ is the $A$-graded module with $V(\alpha)_{\alpha'}
= V_{\alpha+\alpha'}$, then $\HH_0(E-\beta;V(\alpha)) \cong
\HH_0(E-\beta+\alpha;V)(\alpha)$.\qed
\end{lemma}

We shall see that Euler-Koszul homology has the useful property of
detecting ``where'' a module is nonzero, the nonzeroness being
measured in the following sense.

\begin{definition}
Let $V$ be an $A$-graded $\CC[\ttt]$-module.
The set of \emph{true degrees}\/ of $V$ is
\[
  \tdeg(V) = \{\beta \in \ZZ^d: V_\beta \neq 0\}.
\]
The set $\qdeg(V)$ of \emph{quasidegrees}\/ of~$V$ is the Zariski
closure in~$\CC^d$ of its true degrees $\tdeg(V)$.
\end{definition}

Because of the next lemma, we shall often refer to quasidegree sets
as \emph{arrangements}.

\begin{lemma}\label{l:arrangement}
Let $R$ be a noetherian $A$-graded ring that is finitely generated
over its degree~$0$ piece.  The quasidegree set of any finitely
generated graded 
$R$-module is a finite union of affine subspaces
of\/~$\CC^d$, each spanned by the degrees of some subset of the
generators of~$R$.
\end{lemma}
\begin{proof}
Every $A$-graded module has an $A$-graded associated prime, and
therefore a submodule isomorphic to an $A$-graded translate of a
quotient by an $A$-graded prime.  Now use Noetherian induction to
conclude that every such module has a filtration whose successive
quotients are $A$-translates of quotients of $R$ modulo prime ideals.
But being an integral domain, the true degree set of a quotient
$R/\pp$ by a prime ideal~$\pp$ is the affine semigroup generated by
the degrees of the generators of~$R$ that remain nonzero in~$R/\pp$.
\end{proof}

\begin{example}\label{ex:qdeg}
Let $I =
\<bd-de,bc-ce,ab-ae,c^3-ad^2,a^2d^2-de^3,a^2cd-ce^3,a^3d-ae^3\>$ be a
binomial ideal in $\CC[\ttt]$, where we write $\ttt =
(\del_1,\del_2,\del_3,\del_4,\del_5) = (a,b,c,d,e)$, and let
\[
  A = \left[\begin{array}{ccccc}
	1 & 1 & 1 & 1 & 1 \\
	0 & 1 & 2 & 3 & 1
      \end{array}\right]
\qquad\text{and}\qquad
  B = \left[\begin{array}{rrr}
	-2 & -1 &  0 \\ 
	 3 &  0 &  1 \\
	 0 &  3 &  0 \\ 
	-1 & -2 &  0 \\
	 0 &  0 & -1
      \end{array}\right].
\]
One easily verifies that the binomial ideal~$I$ is graded by $\ZZ A =
\ZZ^2$.  If $\omega$ is a primitive cube root of unity ($\omega^3 =
1$), then $I$, which is a radical ideal, has the prime decomposition
\begin{align*}
I =
  \<a,c,d\> &\cap\<bc-ad,b^2-ac,c^2-bd,b-e\>
\\          &\cap\<\omega bc-ad,b^2-\omega ac,\omega^2 c^2-bd,b-e\>
\\          &\cap\<\omega^2bc-ad,b^2-\omega^2ac,\omega c^2-bd,b-e\>.
\end{align*}
If $V = \CC[a,b,c,d,e]/\<a,c,d\>$, then $\qdeg(V)$ is the diagonal
line in~$\CC^2$.  In contrast, the quotient by each one of the other
three prime ideals there has quasidegree set equal to all of~$\CC^2$.
It follows that $\qdeg(\CC[\ttt]/I) = \CC^2$.
\end{example}

Let $\mm$ be the maximal ideal $\<\del_1,\dots \del_n\>$
of~$\CC[\ttt]$.  Since $A$ is pointed with no nonzero columns, $\mm$
is the unique maximal $A$-graded ideal.  Given an $A$-graded
$\CC[\ttt]$-module~$V$, its \emph{local cohomology modules}
\[
  H^i_\mm(V) = \dlim_t \ext^i_{\CC[\ttt]}(\CC[\ttt]/\mm^t,V)
\]
supported at\/~$\mm$ are $A$-graded; see \cite[Chapter~13]{cca}.
Even when $V$ is finitely generated, its local cohomology
modules~$H^i_\mm(V)$ need not be; but their Matlis duals are, so their
quasidegree sets are still arrangements.

\begin{lemma}\label{l:qdeg}
If\/ $V$ is a finitely generated\/ $A$-graded $\CC[\ttt]$-module,
then the quasidegree set $\qdeg(H^i_\mm(V))$ of the $i^\th$ local
cohomology module of\/~$V$ is a union of finitely many integer
translates of the complex subspaces $\CC A_J \subseteq \CC^n$ spanned
by $\{a_j : j \in J\}$ for various~$J$.
\end{lemma}
\begin{proof}
Let $\vea = \sum_{j=1}^n a_j$.  In the graded version
\cite[Theorem~6.3]{gradedGM} of the Greenlees-May theorem \cite{gm},
setting ${\mathcal E}$ equal to the injective hull of the residue
field $\CC[\ttt]/\mm$ yields the natural $A$-graded local duality
vector space isomorphism
\[
  \ext^{n-i}_{\CC[\ttt]}(V,\CC[\ttt])_\alpha \cong
  \hom_\CC(H^i_\mm(V)_{-\alpha+\vea},\CC).
\]
(Use the case ${\mathcal G} = \CC[\ttt]$ to deduce that the right-hand
side of \cite[Theorem~6.3]{gradedGM} is the derived $\hom$ into the
canonical module $\omega_{\CC[\ttt]}$, which is isomorphic as a graded
module to the principal ideal $\<\del_1 \cdots \del_n\>$; see also
\cite[Section~3.5]{BH93} for $\ZZ$-graded local duality.)  Hence $\vea
+ \qdeg(H^i_\mm(V)) = -\qdeg(\ext^{n-i}_{\CC[\ttt]}(V,\CC[\ttt]))$ is
the negative of the quasidegree set of a finitely generated module.
The result is now a consequence of Lemma~\ref{l:arrangement}.
\end{proof}

\section{Binomial primary decomposition}\label{s:primDecomp}

In this section we 
review some prerequisites on primary decomposition of binomial ideals
from \cite{binomialideals} and \cite{primDecomp}, including
interactions with $A$-gradings.  For the applications to Horn
$D$-modules in Section~\ref{s:horn}, we pay special attention to
lattice basis ideals.  For the duration of this section we work
over a polynomial ring $\CC[\ttt]$ in commuting variables~$\ttt =
\del_1,\ldots,\del_n$.

If $L \subseteq \ZZ^n$ is a sublattice, then the \emph{lattice
ideal}\/ of~$L$ is $I_L = \<\ttt^{u_+}-\ttt^{u_-}: u=u_+-u_-
\in L\>$.  Here and henceforth, $u_+$ has $i^{\rm th}$ coordinate
$u_i$ if $u_i\geq 0$ and $0$ otherwise.  The vector $u_- \in \NN^q$ is
defined by $u_+ - u_- = w$, or equivalently, $u_- = (-u)_+$.
More general than $I_L$ are the~ideals
\[
  I_\rho = \<\ttt^{u_+}-\rho(u)\ttt^{u_-}: u=u_+-u_-\in L\>
\]
for any \emph{partial character} $\rho : L \to \CC^*$ of~$\ZZ^n$,
which includes the data of both its domain lattice $L \subseteq \ZZ^n$
and the map to~$\CC^*$.  (The ideal~$I_\rho$ is called $I_+(\rho)$ in
\cite{binomialideals}.)  The ideal $I_\rho$ is prime if and only if
$L$ is a \emph{saturated} sublattice of~$\ZZ^n$, meaning that $L$
equals its \emph{saturation}
\[
  \sat(L) = (\QQ L) \cap \ZZ^n,
\]
where $\QQ L = \QQ \otimes_\ZZ L$ is the rational vector space spanned
by~$L$ in $\QQ^n$.  In fact \cite[Corollary~2.6]{binomialideals},
every binomial prime ideal in $\CC[\ttt]$ has the form
\[
  I_{\rho,J} = I_\rho + \<\del_j : j \notin J\>
\]
for some saturated partial character~$\rho$ (i.e., whose domain is a
saturated sublattice) and subset $J \subseteq \{1,\ldots,n\}$ such
that the binomial generators of~$I_\rho$ only involve variables
$\del_j$ for $j \in J$ (some of which might actually be absent from
the generators of~$I_\rho$).

\begin{example}\label{ex:binomial}
The intersectand $\<a,c,d\>$ in Example~\ref{ex:qdeg} equals the prime
ideal $I_{\rho,J}$ for $J = \{2,5\}$ and $L = \{0\} \subseteq \ZZ^J$.
The remaining three intersectands are the prime ideals $I_{\rho,J}$
for the three characters $\rho$ that are defined on~$\ker(A)$ but
trivial on its index~$3$ sublattice~$\ZZ B$ spanned by the columns~of
$B$, where $J = \{1,2,3,4,5\}$.
\end{example}

\begin{theorem}[{\cite[Theorem~3.2]{primDecomp}}]\label{t:primDecomp}
Fix a binomial ideal $I$.  Write $\ttt_J$ for the monomial
$\prod_{j\in J} \del_j$.  Each associated prime $I_{\rho,J}$ has an
explicitly defined monomial ideal $U_{\rho,J}$ such that
\[
  I = \bigcap_{I_{\rho,J} \in {\rm Ass}(I)} \cC_{\rho,J}
  \quad\text{for}\quad
  \cC_{\rho,J} = \big((I + I_\rho):\ttt_J^\infty\big) + U_{\rho,J}
\]
is a decomposition of $I$ as an intersection of primary binomial ideals.
\end{theorem}
It is not important for our present purposes precisely what
$U_{\rho,J}$ is in general; all we need are various consequences,
especially for the structure of the quotients
$\CC[\ttt]/\cC_{\rho,J}$, derived in \cite{primDecomp} from the
explicit description.  The flavor is captured in the following example
and in Example~\ref{ex:I(B)'}, where the precise answer for certain
minimal primes is quite clean.

\begin{example}\label{ex:I(B)}
Fix matrices $A$ and~$B$ as in Convention~\ref{conv:B}.  This
identifies $\ZZ^d$ with the quotient of $\ZZ^n/\ZZ B$ modulo its
torsion subgroup.  The \emph{lattice basis ideal}\/ corresponding to
the lattice $\ZZ B = \{Bz : z \in \ZZ^m\}$ is defined by
\[
  I(B) = \<\ttt^{u_+}-\ttt^{u_-}: u=u_+-u_- \;\mbox{is a column of}\;
  B\> \subseteq \CC[\del_1,\dots ,\del_n].
\]
Each of the minimal primes of~$I(B)$ arises, after row and column
permutations, from a block decomposition of~$B$ of the form
\begin{equation}\label{eq:MNOB}
  \left[
    \begin{array}{l|r}
    N & B_J\!\\\hline
    M & 0\ 
    \end{array}
  \right],
\end{equation}
where $M$ is a mixed submatrix of~$B$ of size $q \times p$ for some $0
\leq q \leq p \leq m$ \cite{hostenshapiro}.  (Matrices with $q = 0$
rows are automatically mixed; matrices with $q = 1$ row are never
mixed.)  We note that not all such decompositions correspond to
minimal primes: the matrix $M$ has to satisfy another condition which
Ho\c{s}ten and Shapiro call irreducibility \cite[Definition~2.2 and
Theorem~2.5]{hostenshapiro}.  If $I(B)$ is a complete intersection,
then only square matrices $M$ will appear in the block
decompositions~(\ref{eq:MNOB}), by a result of Fischer and Shapiro
\cite{fischer-shapiro}.

For each partial character $\rho : \sat(\ZZ B_J) \to \CC^*$ extending
the trivial character on~$\ZZ B_J$, the ideal $I_{\rho,J}$ is
associated to~$I(B)$, where $J = J(M) = \{1,\ldots,n\} \minus {\rm
rows}(M)$ indexes the $n-q$ rows not in~$M$.  We reiterate that the
symbol $\rho$ here includes the specification of the sublattice
$\sat(\ZZ B_J) \subseteq \ZZ^n$.  The corresponding primary component
$\cC_{\rho,J}$ of~$I(B)$ is simply $I_\rho$ if $q = 0$, but will in
general be non-radical when $q \geq 2$ (recall that $q = 1$ is
impossible).
\end{example}

Since $A$-gradings are central to our theory, we collect some relevant
results from \cite{primDecomp}.  Recall Conventions~\ref{conv:A}
and~\ref{conv:B}.  Henceforth, $A_J$ denotes the submatrix of~$A$
whose columns are indexed by~$J$.  We write $\ZZ A_J \subseteq \ZZ^d =
\ZZ A$ for the group generated by these columns.

\begin{lemma}\label{l:toralAndean}
Fix a partial character $\rho : L \to \CC^*$ for a saturated
sublattice $L \subseteq \ZZ^J \subseteq \ZZ^n$.  Let $\cC_{\rho,J}$ be
an $A$-graded binomial $I_{\rho,J}$-primary ideal.  Then $L \subseteq
\ZZ^J \cap \ker_{\ZZ}(A) = \ker_{\ZZ}(A_J)$, the Krull dimension
satisfies $\dim(\CC[\ttt]/I_{\rho,J}) \geq \rank(A_J)$, and the
following are equivalent.
\begin{itemize}
\item
The \emph{Hilbert function} $\ZZ A \to \NN$ defined by $\alpha \mapsto
\dim_\CC(\CC[\ttt]/\cC_{\rho,J})_\alpha$ is bounded above.
\item
The homomorphism $\ZZ^J/L \onto \ZZ A_J \subseteq \ZZ^d$ is injective.
\item
$L = \ker_\ZZ(A_J)$.
\item
$\dim(\CC[\ttt]/I_{\rho,J}) = \rank(A_J)$.
\end{itemize}
When these conditions are satisfied, the module
$\CC[\ttt]/\cC_{\rho/J}$ and the lattice~$L$ are called \emph{toral},
the ideal $I_{\rho,J}$ is called a \emph{toral prime}, and
$\cC_{\rho,J}$ is called a \emph{toral (primary) component}.  When
these conditions are not satisfied, substitute \emph{Andean}
(see Remark~{\ref{rk:andean}}) for ``toral'' above.
\end{lemma}
\begin{proof}
These conditions are the ones appearing, respectively, in
\cite[Definition~4.3, Proposition~4.7, Corollary~4.8, and
Lemma~4.9]{primDecomp}.
\end{proof}

\begin{example}\label{ex:binomial'}
In Example~\ref{ex:binomial}, the homomorphism $A_{\{2,5\}}: \ZZ^{\{2,5\}} \to \ZZ^2$
is not injective since it maps both basis vectors to
$\left[\twoline{1}{1}\right]$; thus the prime ideal $\<a,c,d\>$ is an
Andean component of~$I$.  In contrast, the remaining associated prime
ideals are all toral by Lemma~\ref{l:toralAndean}, with $A_J = A$.
\end{example}

The final example in this section demonstrates, at long last, just how
concrete binomial primary decomposition can be when expressed in
combinatorial terms.  It will be applied directly in
Section~\ref{s:horn} to construct solutions to binomial Horn systems.
Example~\ref{ex:I(B)'} was, for us, the motivation and starting point
for all of the other results in this article and in \cite{primDecomp}.
To state it, we need a definition.

\begin{definition}\label{d:M}
Any integer matrix $M$ with $q$ rows defines an undirected
graph~$\Gamma(M)$ having vertex set $\NN^q$ and an edge from $u$
to~$v$ if $u - v$ or $v - u$ is a column of~$M$.  An
\emph{$M$-subgraph}\/ of~$\NN^q$ is a connected component
of~$\Gamma(M)$.  An $M$-subgraph is \emph{bounded}\/ if it has
finitely many vertices, and \emph{unbounded}\/ otherwise.  (See
Example~\ref{ex:concrete-subgraph} for an explicit
computation~in~$\NN^3$.)
\end{definition}
  
\begin{example}\label{ex:I(B)'}
Resume the notation of Example~\ref{ex:I(B)}.  If $I_{\rho,J}$ is a
toral minimal prime of~$I(B)$ given by a matrix decomposition as
in~(\ref{eq:MNOB}), so $J = J(M)$, then
\[
  \cC_{\rho,J} = I(B) + I_{\rho,J} + U_M,
\]
where $U_M \subseteq \CC[\del_j : j \in \oJ]$ is the ideal
$\CC$-linearly spanned by all monomials with exponent vectors in the
union of the unbounded $M$-subgraphs of\/~$\NN^\oJ$; this is
\cite[Corollary~4.14]{primDecomp}, which also says that every monomial
in $\cC_{\rho,J}$ already lies in~$U_M$.
\end{example}

\section{Toral modules}\label{s:toral}

Much of this article concerns widely diverging $D$-module theoretic
behavior lifted from the toral vs.\ Andean dichotomy in the primary
components of graded binomial ideals.  The functor translating to
$D$-modules is Euler-Koszul homology, which was originally conceived
of for \emph{toric}\/ modules \cite[Definition~4.5]{MMW}.  Here, we
shall show that all of the main results in~\cite{MMW} hold, with
essentially the same proofs, for the more general class of
\emph{toral} modules in Definition~\ref{d:toral}.  The key starting
point is the filtration characterization in Proposition~\ref{p:toral}.
Our main results for toral modules are Theorems~\ref{t:rigid},
\ref{t:jumpqdeg}, \ref{t:holfamily}, and~\ref{t:rankjumps}.

\begin{definition}\label{d:toral}
An $A$-graded $\CC[\ttt]$-module~$V$ is \emph{natively toral}\/ if
there is a binomial prime ideal~$I_{\rho,J}$ and a degree $\alpha \in
\ZZ^d$ such that $V(\alpha) \cong \CC[\ttt]/I_{\rho,J}$ is a toral
quotient (Lemma~\ref{l:toralAndean}).  The module $V$ is
\emph{toral}\/ if it is finitely generated and its $A$-graded Hilbert
function is bounded.
\end{definition}

\begin{prop}\label{p:toral}
An $A$-graded $\CC[\ttt]$-module~$V$ is toral if and only if it has a
filtration $0 = V_0 \subset V_1 \subset \cdots \subset V_{\ell-1}
\subset V_\ell = V$ whose successive quotients $V_k/V_{k-1}$ are all
natively toral.
\end{prop}
\begin{proof}
The proof proceeds by Noetherian induction to reduce to the prime
case, and then by showing that every toral prime is binomial.  The
argument is the same as for \cite[Proposition~4.7]{primDecomp}, but
with general modules~$V$ in place of primary quotients
$\CC[\ttt]/\cC_{\rho,J}$.%
\end{proof}

The argument in the proof of Proposition~\ref{p:toral} actually shows
more.

\begin{lemma}\label{l:W}
If\/ $W \subseteq V$ are $A$-graded modules with\/ $V$ toral, then
$W$ and\/ $V/W$ are toral.
\end{lemma}
\begin{proof}
Intersecting any toral filtration of~$V$ with $W$ yields a filtration
of~$W$ whose successive quotients are toral because they are
$A$-graded modules over natively toral quotients
$\CC[\ttt]/I_{\rho,J}$.  Hence $W$ is toral.  The same argument works
for the image filtration in~$V/W$.
\end{proof}

We begin recounting the results of \cite{MMW} with an elementary
observation about how Euler-Koszul homology works for modules killed
by some of the variables; the proof is the same as
\cite[Lemma~4.8]{MMW}.  For notation, let $E_i^J$ be the operator
obtained from~$E_i$ by setting the terms $x_j\del_j$ to zero for $j
\notin J$.  This operator can be thought of as lying in the Weyl
algebra $D_J$ in the variables $x_j$ and~$\del_j$ for $j \in J$.
Denote by $x_\oJ$ the $x$-variables for $j \notin J$.

\begin{lemma}\label{l:EJ}
If the variables $\del_j$ for $j \notin J$ annihilate an $A$-graded
$\CC[\ttt]$-module $V$, then $D \otimes_{\CC[\ttt]} V \cong \CC[x_\oJ]
\otimes_\CC (D_J \otimes_{\CC[\ttt_J]} V)$ as $D = D_\oJ \otimes_\CC
D_J$-modules.  Acting by $E_i$ on $D \otimes_{\CC[\ttt]} V$ as
in~(\ref{eq:E}) is the same as acting by $E_i^J$ on the right-hand
factor of\/ $\CC[x_\oJ] \otimes_\CC (D_J \otimes_{\CC[\ttt_J]}
V)$.\qed
\end{lemma}

Many of the following results are stated in the context of
\emph{holonomic} $D$-modules, which by definition are the finitely
generated left $D$-modules~$W$ with \mbox{$\ext^j_D(W,D) = 0$}
for~\mbox{$j \neq n$}.  When $W$ is holonomic, the vector
space~$\CC(x) \otimes_{\CC[x]} W$ over the field $\CC(x)$ of rational
functions in $x_1,\ldots,x_n$ has finite dimension equal to the
\emph{holonomic rank}\/~$\rank(W)$ by a celebrated theorem of
Kashiwara; see \cite[Theorem~1.4.19 and Corollary~1.4.14]{SST}.

We shall also be interested in whether our $D$-modules are
\emph{regular holonomic}, the definition of which can be found in
\cite{Bjork}.  For an $A$-hypergeometric $D$-module
(Example~\ref{ex:A-hyp}), regular holonomicity is known
\cite{equivariant} to occur when $A$ is \emph{homogeneous}, meaning
that there is a row vector $\psi \in \QQ^d$ such that $\psi A$ equals
the row vector $[1,\ldots,1]$.  In this case, the $\ZZ A=\ZZ^d$-grading on
$\CC[\ttt]$ coarsens naturally to the \emph{standard $\ZZ$-grading},
in which $\deg(\del_j) = 1 \in \ZZ$ for all~$j$.

\begin{theorem}\label{t:rigid}
If $V$ is a toral $\CC[\ttt]$-module and $\beta \in \CC^d$, then the
Euler-Koszul homology $\HH_i(E-\beta;V)$ is holonomic for all~$i$.
Moreover, the following are equivalent.
\begin{numbered}
\item
$\HH_0(E-\beta;V)$ has holonomic rank~$0$.

\item
$\HH_0(E-\beta;V) = 0$.

\item
$\HH_i(E-\beta;V) = 0$ for all $i \geq 0$.

\item
$-\beta \not\in \qdeg(V)$.
\end{numbered}
If, in addition, the matrix~$A$ is homogeneous, then
$\HH_i(E-\beta;V)$ is regular holonomic for~all~$i$.
\end{theorem}
\begin{proof}
This is the toral generalization of \cite[Proposition~5.1]{MMW} and
\cite[Proposition~5.3]{MMW}.  To see that it holds, start with
\cite[Notation~4.4]{MMW}: instead of only allowing submatrices of~$A$
corresponding to faces of the semigroup~$\NN A$, we allow
submatrices~$A_J$ with arbitrary column sets $J \subseteq
\{1,\ldots,n\}$.  Then, in \cite[Definition~4.5]{MMW}, replace
``toric'' with ``toral'' and change $S_{F_k}$ to
$\CC[\ttt]/I_{\rho,J}$; that this defines toral modules is by
Lemma~\ref{l:toralAndean}.

The key is \cite[Lemma 4.9]{MMW}.  In the proof there, first replace
${\mathcal M}^F_\beta = D/H_A(I^F_A,\beta)$ by
$D/H_A(I_{\rho,J},\beta)$.  Then observe that rescaling the variables
via~$\rho$ induces an $A$-graded automorphism of~$D$ commuting with
the construction of Euler-Koszul complexes (because $x_j\del_j$ is
invariant under the automorphism).  Hence the theorem for natively
toral modules need only be proved in the special case $\rho =
\mbox{}$identity.  This allows us to use $I_{A_J} + \<\del_j:j\not
\in J\>$ instead of~$I_{\rho,J}$.  The rest of the proof of
\cite[Lemma 4.9]{MMW} goes through unchanged, and when $A$ is
homogeneous, provides regular holonomicity as a consequence of the
analogous result for GKZ systems from \cite{equivariant,uli06}.

Now extend the proof of \cite[Proposition 5.1]{MMW} to the toral
setting.  For the first paragraph of that proof, replace ``toric''
with ``toral'' and replace $S_F$ by $\CC[\ttt]/I_{\rho,J}$.  For the
later paragraphs of the proof, begin by working with the module~$M$
there being native toral.  This allows us to replace~$I_A$, when it
arises as an annihilator toward the end, with $I_{\rho,J}$, thereby
proving the native toral case.  For the arbitrary toral case, simply
note that for any exact sequence $0 \rightarrow V' \rightarrow V
\rightarrow V'' \rightarrow 0$ in which $V'$ and $V''$ both have
(regular) holonomic Euler-Koszul homology, each Euler-Koszul homology
module of $V$ is placed between two (regular) holonomic modules, and
is hence (regular) holonomic.

Finally, to generalize \cite[Proposition~5.3]{MMW}, replace ``toric''
with ``toral'' in the statement and proof.  Then, in the proof,
replace $I^F_A$ by $I_{\rho,J}$ and $S_F$ by $\CC[\ttt]/I_{\rho,J}$.
\end{proof}

Next we record the toral generalization of \cite[Theorem~6.6]{MMW}.

\begin{theorem}\label{t:jumpqdeg}
The Euler-Koszul homology $\HH_i(E-\beta;V)$ of a toral module~$V$ is
nonzero for some $i \geq 1$ if and only if $-\beta \in
\qdeg(H^i_\mm(V))$ for some $i < d$.  More precisely, if $k$ equals
the smallest homological degree~$i$ for which $-\beta \in
\qdeg(H^i_\mm(V))$, then $\HH_{d-k}(E-\beta;V)$ is holonomic of
nonzero rank while $\HH_i(E-\beta;V)=0$ for $i>d-k$.
\end{theorem}
\begin{proof}
Begin by noting that $\ext^i_{\CC[\ttt]}(V,\CC[\ttt])$ is toral
whenever $V$ is toral.  This is the toral generalization of
\cite[Lemma~6.1]{MMW}; the same proof works, mutatis mutandis,
replacing $S_A$ in \cite{MMW} by $\CC[\ttt]/I_{\rho,J}$ here.  Now
extend \cite[Theorem~6.3]{MMW} to the toral case: the only property of
toric modules used in its proof is the holonomicity of Euler-Koszul
homology, which we have shown is true for toral modules in
Theorem~\ref{t:rigid}.  Finally, to torally extend the toric
\cite[Theorem~6.6]{MMW}, start with the first sentence of the proof,
which for toral modules is Lemma~\ref{l:dim}, below.  After that, the
proof goes through verbatim, given that we have shown the results it
cites for toric modules to be true for toral modules.
\end{proof}

\begin{lemma}\label{l:dim}
If\/ $V$ is toral, then its Krull dimension satisfies $\dim(V) =
\dim(\qdeg(V)) \leq d$.
\end{lemma}
\begin{proof}
For natively toral modules this follows from
Lemma~\ref{l:toralAndean}.  For arbitrary toral modules, the Krull
dimension and the dimension of the quasidegree set both equal the
maximum of the corresponding dimensions for the composition factors in
any toral filtration.
\end{proof}

One of the observations in \cite{MMW} is that hypergeometric systems
$D/H_A(I,\beta)$ for varying~$\beta$ should be viewed as a family of
$D$-modules fibered over~$\CC^d$.  If (the holonomic rank function of
the $D$-modules in) such a family is to behave well, it suffices to
verify that it is a \emph{holonomic family}
\cite[Definition~2.1]{MMW}.  For families arising from toric modules
this is done in \cite[Theorem~7.5]{MMW}, which we now generalize to
the toral setting.  As a matter of notation, let $b = b_1,\ldots,b_d$
be commuting variables of degree zero, so $D[b]$ is a polynomial
algebra over the Weyl algebra~$D$.  For any $A$-graded
$\CC[\ttt]$-module~$V$, construct the \emph{global Euler-Koszul
complex} $\KK_\spot(E-b;V)$ of left $D[b]$-modules and \emph{global
Euler-Koszul homology} $\HH_\spot(E-b;V)$ by replacing $D$ and
$\beta$ in Definition~\ref{d:EK} with $D[b]$ and~$b$ here.  Finally,
if $\CC(x)$ is the field of rational functions in $x_1,\ldots,x_n$,
write $\VV(x) = \CC(x) \otimes_{\CC[x]} \VV$ for any
$\CC[x]$-module~$\VV$, including $\VV = \CC[b][x]$, where we set
$\VV(x) = \CC[b](x)$.

\begin{theorem}\label{t:holfamily}
If\/ $V$ is toral, then the sheaf $\tilde{\VV}$ on~$\CC^d$ whose
global section module is $\VV=\HH_0(E-b;V)$ constitutes a holonomic
family over\/ $\CC^d$; in other words, $\VV_\beta = \HH_0(E-\beta;V)$
is holonomic for all $\beta \in \CC^d$, and\/ $\VV(x)$ is finitely
generated as a module over $\CC[b](x)$.
\end{theorem}
\begin{proof}
\cite[Proposition~7.4]{MMW} holds for $I_{\rho,J}$ in place of $I^F_A$
after harmlessly rescaling the $x$ and~$\xi$ variables inversely to
each other, which affects neither $Ax\xi$ nor the initial ideal in
question.  Therefore we may, in the proof of \cite[Theorem~7.5]{MMW},
simply change ``toric'' to ``toral'' and base the induction again on
$I_{\rho,J}$ and $\CC[\ttt]/I_{\rho,J}$ instead of $I^F_A$ and~$S_A$.
\end{proof}

Considering $b_i$ and $\beta_i$ as elements in the polynomial
ring~$\CC[b]$, we can take ordinary Koszul homology
$H_\spot(b-\beta;W)$ for any $\CC[b]$-module~$W$.  This gets used in
the generalization of \cite[Theorem~8.2]{MMW} to arbitrary
$A$-graded $\CC[\ttt]$-modules, which we state along with the
toral generalization of \cite[Theorem~9.1]{MMW}.  For the latter, we
need also the \emph{jump arrangement} $\cZ_\jump(V) = \bigcup_{i \leq
d-1} \qdeg(H^i_\mm(V))$ of an $A$-graded module~$V$
over~$\CC[\ttt]$.

\begin{theorem}\label{t:rankjumps}
If\/ $V$ is an $A$-graded\/ $\CC[\ttt]$-module and\/ $\VV =
\HH_0(E-b;V)$, then
\[
  \HH_i(E-\beta;V) \cong H_i(b-\beta;\VV),
\]
the left and right sides being Euler-Koszul and ordinary Koszul
homology, respectively.  If,~in addition, $V$ is toral, then $-\beta$
lies in the jump arrangement $\cZ_\jump(V)$ if and only if the
holonomic rank of\/ $\HH_0(E-\beta;V)$ is not minimal (among all
possible choices of~$\beta$).
\end{theorem}
\begin{proof}
\cite[Theorem~8.2]{MMW} and its proof both work verbatim for arbitrary
$A$-graded $\CC[\ttt]$-modules.  That being given, the proof of
\cite[Theorem~9.1]{MMW} works just as well for toral modules, since we
have now seen that all of the earlier results in \cite{MMW} do.
\end{proof}

\section{Andean modules}\label{s:andean}

The finiteness properties of toral modules encapsulated by
Theorem~\ref{t:rigid} will be contrasted in Corollary~\ref{c:rigid}
(the heart of which is Theorem~\ref{t:andeanBig}) with the
infiniteness that occurs for Andean modules.  The feature of toral
modules that drives the proofs in Section~\ref{s:toral} is the toral
filtration in Proposition~\ref{p:toral}.  It would be optimal if we
could simply define an Andean module, in general, to mean one that is
not toral---that is, one whose Hilbert function is unbounded---and
conclude a similar filtration feature for Andean modules.  Alas, this
notion of Andean module is too inclusive for our purposes: it does not
imply a filtration characterization, in general, even though for the
quotient of~$\CC[\ttt]$ by a binomial primary ideal, the unbounded
Hilbert function characterization is equivalent to the filtration one
(Example~\ref{ex:andean}).  Therefore, we take as our foundation the
filtration feature.  The particular form of this feature is dictated
by combinatorial primary decomposition, particularly
\cite[Example~4.6]{primDecomp}.

\begin{definition}\label{d:andean}
An $A$-graded $\CC[\ttt]$-module~$V$ is \emph{natively Andean}\/ if
there is an $\alpha \in \ZZ^d$ and an Andean quotient ring
$\CC[\ttt]/I_{\rho,J}$ (Lemma~\ref{l:toralAndean}) over which
$V(\alpha)$ is torsion-free of rank~$1$ and admits a $\ZZ^J/L$-grading
that refines the $A$-grading via $\ZZ^J/L \to \ZZ^d=\ZZ A$, where
$\rho$ is defined on~$L \subseteq \ZZ^J$.  If $V$ has a finite
filtration $0 = V_0 \subset V_1 \subset \cdots \subset V_{\ell-1}
\subset V_\ell = V$ whose successive quotients $V_k/V_{k-1}$ are all
natively Andean, then $V$ is \emph{Andean}.
\end{definition}

\begin{example}\label{ex:andean}
$\CC[\ttt]/\cC_{\rho,J}$ is Andean for any Andean primary component
$\cC_{\rho,J}$ of any $A$-graded binomial ideal.  This follows
immediately from the statements
\cite[Corollaries~2.13 and~3.3]{primDecomp} about gradings and
filtrations for primary binomial ideals.
\end{example}

\begin{remark}\label{rk:andean}
The adjective ``Andean'' describes the geometry of the gradings on the
$\CC[\ttt]$-modules $\CC[\ttt]/\cC_{\rho,J}$: collapsing (coarsening)
the natural grading by the $\ZZ^J/L$-torsor $\cB$ to the $A$-grading
\cite[Corollary~2.13]{primDecomp} makes the $\cB$-graded degrees sit
like a high thin mountain range over~$\ZZ^d$, supported on finitely
many translates of~$\ZZ A_J$.
\end{remark}

Here is a weak form of Euler-Koszul rigidity for Andean modules (but
see Corollary~\ref{c:rigid}).

\begin{lemma}\label{l:weak-rigidity}
If\/ $V$ is an Andean module and $-\beta \not \in \qdeg(V)$, then
$\HH_i(E-\beta;W) = 0 = \HH_i(E-\beta;V/W)$\/ for all~$i$ and all
$A$-graded submodules $W \subseteq V$.
\end{lemma}
\begin{proof}
First assume that $V$ is natively Andean.  The torsion-freeness
ensures that $\qdeg(V)$ is a $\ZZ^d$-translate of the complex span
$\CC A_J$ of the columns of~$A$ indexed by~$J$, so let us also assume
for the moment that \mbox{$\qdeg(V) = \CC A_J$}.  The result for
this~$V$ and all of its $A$-graded submodules follows from
Lemma~\ref{l:EJ}, because the $\CC$-linear span of $E_1^J-\beta_1,
\ldots, E_d^J-\beta_d$ contains a nonzero scalar if $\beta \notin \CC
A_J$ (some linear combination of $E_1^J,\ldots,E_d^J$ is zero, while
the corresponding linear combination of $\beta_1,\ldots,\beta_d$ is
nonzero, and hence a unit).

The case where $V$ is natively Andean (or a submodule thereof) and
$\qdeg(V) = \alpha + \CC A_J$ is proved by applying the above argument
to~$V(-\alpha)$, using Lemma~\ref{l:V}.  The case where $V$ is a
general Andean module is proved by induction on the length of an
Andean filtration, using that $\qdeg(V) = \qdeg(V') \cup \qdeg(V'')$
whenever $0 \to V' \to V \to V'' \to 0$ is an exact sequence.
Finally, for an $A$-graded submodule~$W$ of a general Andean
module~$V$, intersecting~$W$ with an Andean filtration of~$V$ yields a
filtration of~$W$ whose successive quotients are submodules of native
Andean modules.  Hence the proof of vanishing of Euler-Koszul homology
by induction on the length of the filtration still applies.

The vanishing of all $\HH_i(E-\beta;V/W)$ follows easily from the
vanishing for~$V$ and for~$W$.
\end{proof}

The following lemma will allow us to reduce to the case of Andean
quotients $\CC[\ttt]/I_{\rho,J}$ whenever we need to work with
natively Andean modules.

\begin{lemma}\label{l:andeanToral}
A natively Andean module $V$ has a filtration whose successive
quotients are $A$-graded translates of various quotients
$\CC[\ttt]/I_{\rho,J}$, each being natively either toral or Andean. At
least one of these quotients is natively Andean.
\end{lemma}
\begin{proof}
By definition, $V$ is torsion free of rank~$1$ over an Andean quotient
$\CC[\ttt]/I_{\tau,S}$ where $S \subseteq \{ 1,\dots,n\}$ and $\tau$
is a partial character (we use non-standard notation to avoid
confusion with the statement we need to prove).  Harmlessly rescaling
the variables, we may assume that $I_{\tau,S} = I_L + \<\del_j : j
\notin S\>$, so $\CC[\ttt]/I_{\tau,S}$ is a semigroup ring~$\CC[Q]$
for some $Q \subseteq \ZZ^{S\!}/L$.  Replacing $V$ with an $A$-graded
translate, we may further assume that $V$ is $\ZZ^{S\!}/L$-graded.
Using Noetherian induction as in the proof of
Lemma~\ref{l:arrangement}, we construct a filtration of~$V$ whose
successive quotients are $\ZZ^{S\!}/L$-graded translates of quotients
$\CC[Q]/\pp_{Q'}$ modulo prime ideals~$\pp_{Q'}$ for faces~$Q'
\subseteq Q$ (these are the $\CC[\ttt]/I_{\rho,J}$ of the statement).
Each of these, being $A$-graded, is either natively toral or natively
Andean.  Moreover, if all of them were toral, then $V$ would be toral
as well, so the last assertion follows.
\end{proof}

Now we combine the Euler-Koszul theory for Andean and toral modules to
conclude that the hypergeometric $D$-modules associated to Andean
modules, if nonzero, are very large.

\begin{theorem}\label{t:andeanBig}
If an $A$-graded $\CC[\ttt]$-module $V$ possesses a surjection to
an Andean module~$W$, and if $-\beta \in \qdeg(W)$, then
$\HH_0(E-\beta;V)$ has uncountably many linearly independent solutions
near any general point $x \in \CC^n$; that is, ${\rm
Hom}_D(\HH_0(E-\beta;V),\OO_x)$ is a vector space of uncountable
dimension over\/~$\CC$, where $\OO_x$ is the local ring of analytic
germs~at~$x$.
\end{theorem}

\begin{proof}
Since a surjection of $\CC[\ttt]$-modules induces a surjection of
zeroth Euler-Koszul homology~$\HH_0(E-\beta;\blank)$, we may assume
that $V = W$ is Andean.

Consider an exact sequence $0 \to V' \to V \to V'' \to 0$ of Andean
modules in which $V''$ is natively Andean.  If $-\beta \notin
\qdeg(V'')$, then $\HH_0(E-\beta;V') \cong \HH_0(E-\beta;V)$ by
Lemma~\ref{l:weak-rigidity} for~$V''$, so we may harmlessly replace
$V$ with~$V'$.  Continuing in this manner, using induction on the
length of an Andean filtration of~$V$, we may assume that $-\beta \in
\qdeg(V'')$.  But then, since $\HH_0(E-\beta;V)$ always surjects onto
$\HH_0(E-\beta;V'')$, we may assume that $V = V''$ is natively Andean.
By Lemma~\ref{l:V}, we may further assume that $V$ is torsion-free of
rank~$1$ over some Andean quotient $R = \CC[\ttt]/I_{\rho,J}$, and
that $V$ contains $R$ with no $A$-graded translation.

Using Lemma~\ref{l:andeanToral} and its notation, take a filtration $0
= V_0 \subset V_1 \subset \cdots \subset V_\ell = V$ in which each of
the successive quotients $V_k/V_{k-1}$ is a $A$-graded translate
of some prime quotient that is natively either toral or Andean.  We
are free to choose $V_1 = R$, and we do so.  Let~$k$ be the largest
index such that $V_k/V_{k-1}$ is Andean and $-\beta \in
\qdeg(V_k/V_{k-1})$, noting that such an index exists because $V_1/V_0
= R$ satisfies the condition.  Since $V$ surjects onto $V/V_{k-1}$, we
find that $\HH_0(E-\beta;V)$ surjects onto $\HH_0(E-\beta;V/V_{k-1})$.
Therefore, replacing $V$ by $V/V_{k-1}$ and $\ZZ^d$-translating again
via Lemma~\ref{l:V} if necessary, it is enough to prove the case
$k=1$, with $V_1 = R$.

If the above filtration has length $\ell > 1$, then the kernel and
cokernel of the homomorphism $\HH_0(E-\beta;V_{\ell-1}) \to
\HH_0(E-\beta;V)$ are holonomic, being $\HH_i(E-\beta;V/V_{\ell-1})$
for $i \in \{0,1\}$; this is by Theorem~\ref{t:rigid} if
$V/V_{\ell-1}$ is toral, and by Lemma~\ref{l:weak-rigidity} if
$V/V_{\ell-1}$ is Andean with $-\beta \notin \qdeg(V/V_{\ell-1})$.
Therefore the desired result holds for~$V$ if and only if it holds
for~$V_{\ell-1}$.  This argument reduces us to the case $\ell=1$ by
induction on~$\ell$, so we may assume that $V = R =
\CC[\ttt]/I_{\rho,J}$.

The condition $-\beta \in \qdeg(R)$ means exactly that $-\beta$, or
equivalently~$\beta$, lies in the complex column span~$\CC A_J$.  Let
$\hat A$ be a matrix for the projection $\ZZ^J \to \ZZ^J/L$, and write
$\ZZ\hat A = \ZZ^J/L$.  If $\hat\beta$ is a vector in~$\CC\hat A$
mapping to~$\beta$ under the surjection to~$\CC A_J$ afforded by
Lemma~\ref{l:toralAndean}, then denote by $\hat E - \hat\beta$ the
sequence of Euler operators associated to~$\hat A$.  Thought of as
elements in the space of affine linear functions $\ZZ^J \to \CC$, the
Euler operators $E_1^J-\beta_1,\ldots,E_d^J-\beta_d$ truncated from
$E-\beta$ generate a sublattice $\ZZ\{E^J-\beta\}$ properly contained
in the sublattice $\ZZ\{\hat E-\hat\beta\}$ generated by $\hat E -
\hat\beta$.  The binomial hypergeometric system $D/H_{\hat
A}(I_{\rho,J},\hat\beta)$ is holonomic of positive rank by
Theorem~\ref{t:rigid} (for $\ZZ^J/L$-graded toral
$\CC[\ttt_J]$-modules, via Lemma~\ref{l:EJ}).  Its solutions are also
solutions of $D/H_A(I_{\rho,J},\beta)$ because
\[
  H_A(I_{\rho,J},\beta)
  =
  D \cdot \<I_{\rho,J},\ZZ\{E^J-\beta\}\>
  \subseteq
  D \cdot \<I_{\rho,J},\ZZ\{\hat E-\hat\beta\}\>
  =
  H_{\hat A}(I_{\rho,J},\hat\beta).
\]
On the other hand, for any pair of distinct lifts $\hat\beta \neq
\hat\beta'$, the linear span of $\ZZ\{\hat E-\hat\beta\}$ together
with $\ZZ\{\hat E-\hat\beta'\}$ contains a nonzero scalar.  It follows
that the solutions to $D/H_{\hat A}(I_{\rho,J},\hat\beta)$ for
varying~$\hat\beta$ are linearly independent.  The direct sum of these
(local) solution spaces is therefore an uncountable-dimensional subspace of
the (local) solutions to $\HH_0(E-\beta;R) = D/H_A(I_{\rho,J},\beta)$.
\end{proof}

Summarizing the above results, let us emphasize the dichotomy between
toral and Andean modules by recording the Andean analogue of
Theorem~\ref{t:rigid}.

\begin{corollary}\label{c:rigid}
The following are equivalent for an Andean $\CC[\ttt]$-module $V$ and
$\beta \in \CC^d$.
\begin{numbered}\addtocounter{enumi}{-1}
\item
$\HH_0(E-\beta;V)$ has countable-dimensional local solution space.

\item
$\HH_0(E-\beta;V)$ has finite-dimensional local solution space.

\item
$\HH_0(E-\beta;V) = 0$.

\item
$\HH_i(E-\beta;V) = 0$ for all $i \geq 0$.

\item
$-\beta \not\in \qdeg(V)$. \qed
\end{numbered}
\end{corollary}

\section{Binomial $D$-modules}\label{s:binomial}

Using the functoriality of Euler-Koszul homology, we now deduce the
holonomicity, regularity, and other structural properties of arbitrary
binomial $D$-modules, including the binomial Horn systems which
motivated and presaged the developments here.  Our first principal
result is the specification, for any $A$-graded binomial
ideal~$I$, of an arrangement of finitely many affine subspaces
of~$\CC^d$ such that the binomial $D$-module $D/H_A(I,\beta)$ is
holonomic precisely when $-\beta$ lies outside of it
(Theorem~\ref{t:holonomic}).  Moreover, holonomicity occurs if and
only if the vector space of local solutions to $H_A(I,\beta)$ has finite
dimension.  The subspace arrangement arises from the primary
decomposition of~$I$ into its toral and Andean components.  When
$D/H_A(I,\beta)$ is holonomic,
it is also regular
holonomic if and only if
$I$ is $\ZZ$-graded in the standard sense.
Finally, we construct another finite affine subspace arrangement
in~$\CC^d$ such that for $-\beta$ outside of it, the binomial
$D$-module splits as a direct sum of primary toral binomial
$D$-modules (Theorem~\ref{t:sum}).

For the duration of this section, fix an $A$-graded binomial ideal $I
\subseteq \CC[\del_1,\ldots,\del_n]$ and fix an irredundant primary
decomposition as in Theorem~\ref{t:primDecomp}.  Thus, as in
Lemma~\ref{l:toralAndean}, some of the quotients
$\CC[\ttt]/\cC_{\rho,J}$ are toral and some are Andean.  Much of what
we do is independent of the particular primary decomposition, since
the data we typically need come from the quasidegrees of certain
related modules.  For example, the holonomicity in
Theorem~\ref{t:holonomic} is clearly independent of the primary
decomposition.

\begin{definition}\label{d:Z}
The \emph{Andean arrangement}\/ $\cZ_\andean(I)$ is the
union of the quasidegree sets $\qdeg(\CC[\ttt]/\cC_{\rho,J})$ for the
Andean primary components~$\cC_{\rho,J}$ of~$I$.
\end{definition}

\begin{lemma}\label{l:Z}
The Andean arrangement $\cZ_\andean(I)$ is a union of
finitely many integer translates of the subspaces $\CC A_J \subseteq
\CC^n$ for which there is an Andean associated prime $I_{\rho,J}$.
\end{lemma}
\begin{proof}
Apply Lemma~\ref{l:arrangement} to an Andean filtration of each Andean
component $\CC[\ttt]/\cC_{\rho,J}$.
\end{proof}

\begin{theorem}\label{t:holonomic}
Given the $A$-graded binomial ideal $I \subseteq \CC[\ttt]$, the
following are equivalent.
\begin{numbered}\addtocounter{enumi}{-1}
\item
The vector space of local solutions to $H_A(I,\beta)$ has countable
dimension.
\item
The vector space of local solutions to $H_A(I,\beta)$ has finite dimension.
\item
The binomial $D$-module $D/H_A(I,\beta)$ is holonomic.
\item
The Euler-Koszul homology $\HH_i(E-\beta;\CC[\ttt]/I)$ is holonomic
for all~$i$.
\item
$-\beta \not\in \cZ_\andean(I)$.
\end{numbered}
For $I$ standard $\ZZ$-graded, these are equivalent to regular
holonomicity of\/ $\HH_i(E-\beta;\CC[\ttt]/I)$.
Moreover, the existence of a parameter $\beta$ for which
$\HH_0(E-\beta;\CC[\ttt]/I)$ is regular holonomic
implies that $I$ is $\ZZ$-graded.
\end{theorem}
\begin{proof}
The last claim follows from the rest by Theorem~\ref{t:rigid}
and results in~\cite{equivariant,uli06}.  Item~1
trivially implies item~0.  Item~2 implies item~1 because holonomic
systems have finite rank.  Item~3 implies item~2 by
Definition~\ref{d:Eulers} and Example~\ref{ex:A-hyp}.  If $-\beta \in
\cZ_\andean(I)$, then $-\beta \in \qdeg(\CC[\ttt]/\cC_{\rho,J})$ for
some Andean component $\cC_{\rho,J}$, so item~0 implies item~4 by
Theorem~\ref{t:andeanBig} for the surjection $\CC[\ttt]/I \onto
\CC[\ttt]/\cC_{\rho,J}$.  Finally, item~4 implies item~3 by
Theorem~\ref{t:rigid} and Proposition~\ref{p:Itoral}, below, given
that $\CC[\ttt]/\bigcap_{I_{\rho,J}\ \toral} \cC_{\rho,J}$ 
is a submodule of $\bigoplus_{I_{\rho,J}\
\toral} \CC[\ttt]/\cC_{\rho,J}$ and is hence toral.%
\end{proof}

\begin{prop}\label{p:Itoral}
Let $I_\toral = \bigcap_{I_{\rho,J}\ \toral} \cC_{\rho,J}$ be the
intersection of the toral primary components of~$I$.  If $-\beta$ lies
outside of the Andean arrangement of~$I$, then the natural surjection
$\CC[\ttt]/I \onto \CC[\ttt]/I_\toral$ induces an isomorphism in
Euler-Koszul homology:
\[
  \HH_i(E-\beta;\CC[\ttt]/I) \cong \HH_i(E-\beta;\CC[\ttt]/I_\toral)
  \text{\rm\ for all }i \text{\rm\ when } \hbox{$-$}\beta \notin
  \cZ_\andean(I).
\]
\end{prop}
\begin{proof}
If $I_\andean$ is the intersection of the Andean primary components
of~$I$, then
\[
  \frac{I_\toral}{I}\ =\ \frac{I_\toral}{I_\toral \cap I_\andean}\
  \cong\ \frac{I_\toral + I_\andean}{I_\andean}
\]
is a submodule of $\CC[\ttt]/I_\andean$, which in turn is a submodule
of $\bigoplus_{I_{\rho,J}\ \andean} \CC[\ttt]/\cC_{\rho,J}$.  Since
$\cZ_\andean(I)$ is the quasidegree set of this Andean direct sum, the
exact sequence
\[
  0 \to \frac{I_\toral}{I} \to \frac{\CC[\ttt]}{I} \to
  \frac{\CC[\ttt]}{I_\toral} \to 0
\]
yields isomorphisms $\HH_i(E-\beta;\CC[\ttt]/I) \cong
\HH_i(E-\beta;\CC[\ttt]/I_\toral)$ of Euler-Koszul homology for
all~$i$, by Lemma~\ref{l:weak-rigidity} for $I_\toral/I$.
\end{proof}

Now we move on to the question of when $D/H_A(I,\beta)$ splits into a
direct sum.

\begin{definition}\label{d:P}
The \emph{primary cokernel module}\/ $P_I$ is defined by the exact
sequence
\[
  0 \to \frac{\CC[\ttt]}{I} \to \bigoplus_{I_{\rho,J} \in {\rm Ass}(I)}
  \frac{\CC[\ttt]}{\cC_{\rho,J}} \to P_I \to 0.
\]
The \emph{primary arrangement}\/ is $\cZ_\primary(I) =
\qdeg(P_I) \cup \cZ_\andean(I)$.
\end{definition}

\begin{prop}\label{p:Z}
The primary arrangement $\cZ_\primary(I)$ is a union of finitely many
integer translates of subspaces $\CC A_J \subseteq \CC^n$.  If there
exists $\beta \in \CC^d$ such that the local solution space
of~$H_A(I,\beta)$ has finite dimension, then $\cZ_\primary(I)$ is a
proper Zariski-closed subset of\/~$\CC^d$.
\end{prop}
\begin{proof}
The first sentence is by Lemma~\ref{l:arrangement}.  For the second
sentence, let $(P_I)_\toral$ be the image in~$P_I$ of the direct sum
$\bigoplus_\toral \CC[\ttt]/\cC_{\rho,J}$.  A point in $\qdeg(P_I)$
that does not lie in~$\cZ_\andean(I)$ must necessarily be a
quasidegree of $(P_I)_\toral$; that is
\begin{equation}\label{eq:Z}
  \cZ_\primary(I)=\qdeg\big((P_I)_\toral\big)\cup\cZ_\andean(I).
\end{equation}
The existence of our $\beta$ immediately implies that $\cZ_\andean(I)$
is a proper Zariski-closed subset of~$\CC^n$, so by~(\ref{eq:Z}) we
need only prove the same thing for $\qdeg((P_I)_\toral)$.  The module
$(P_I)_\toral$ is supported on the union of the toric subvarieties
$T_{\rho,J} = {\rm Spec}( \CC[\ttt]/I_{\rho,J})$ for the toral
associated primes of~$I$; this much is by definition.  However, the
map $\CC[\ttt]/I \to \bigoplus_{{\rm Ass}(I)} \CC[\ttt]/\cC_{\rho,J}$
is an isomorphism locally at a point~$x$ whenever $x$ lies in only one
of the associated varieties~$T_{\rho,J}$ (toral or otherwise).
Therefore $(P_I)_\toral$ is supported on the union of the pairwise
intersections of the toral toric varieties~$T_{\rho,J}$ associated
to~$I$.  Hence it is enough to show that if $R$ is the coordinate ring
of the intersection $T_{\rho,J} \cap T_{\rho',J'}$ of any two distinct
toral varieties, then $\qdeg(R)$ is a proper Zariski-closed subset
of~$\CC^d$.  This is a consequence of Lemma~\ref{l:dim}.
\end{proof}

\begin{example}\label{ex:primary}
In Examples~\ref{ex:qdeg}, \ref{ex:binomial}, and~\ref{ex:binomial'},
the primary arrangement $\cZ_\primary(I)$ consists of the five bold
lines in Figure~\ref{f:primary}.  The diagonal line through
$-\!\left[\twoline 11\right]$ is the Andean
arrangement~$\cZ_\andean(I)$ by Examples~\ref{ex:qdeg}
and~\ref{ex:binomial}.
\begin{figure}
\includegraphics{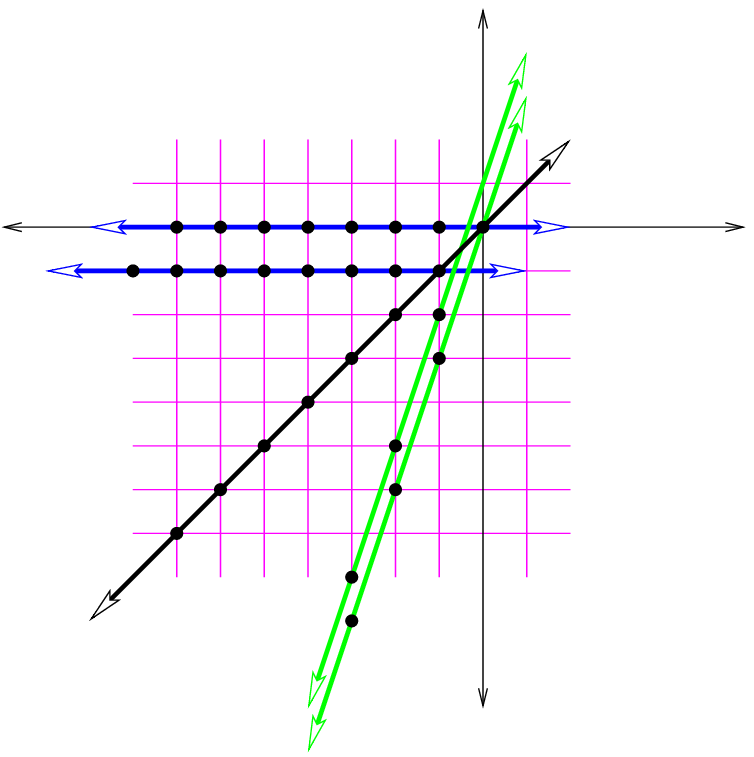}
\caption{Primary arrangement $\cZ_\primary(I)$ of the binomial
ideal~$I$ in Example~\ref{ex:primary}}
\label{f:primary}
\end{figure}
On the other hand, the pairwise intersections of the toral components
of~$I$ all equal $\<bc,ad,b^2,ac,c^2,bd,b-e\>$, which has primary
\mbox{decomposition}
\[
  \<bc,ad,b^2,ac,c^2,bd,b-e\> = \<b^2,c,d,b-e\> \cap \<a,b,c^2,b-e\>.
\]
The set of true degrees of~$P_I$ that lie outside of $\cZ_\andean(I)$
coincides with the true degree set
$\tdeg(\CC[a,b,c,d,e]/\<bc,ad,b^2,ac,c^2,bd,b-e\>)$, which consists
simply of the $A$-degrees of the monomials in $a$, $b$, $c$, and~$d$
that are nonzero in this quotient.  The exponent vectors of these
monomials are those of the form
\[
  \left[\begin{array}{c} \alpha \\ 0 \\ 0 \\ 0 \end{array}\right],
  \left[\begin{array}{c} \alpha \\ 1 \\ 0 \\ 0 \end{array}\right],
  \left[\begin{array}{c} 0 \\ 0 \\ 0 \\ \delta \end{array}\right],
  \text{ or }
  \left[\begin{array}{c} 0 \\ 0 \\ 1 \\ \delta \end{array}\right]
\]
for $\alpha \in \NN$ and $\delta \in \NN$, so $\tdeg(P_I) \minus
\cZ_\andean(I)$ consists of the lattice points having the form
\[
  \left[\begin{array}{c} -\alpha  \\  0        \end{array}\right],
  \left[\begin{array}{c} -\alpha-1\\ -1        \end{array}\right],
  \left[\begin{array}{c} -\delta  \\ -3\delta  \end{array}\right],
  \text{ or }
  \left[\begin{array}{c} -\delta-1\\-3\delta-2 \end{array}\right],
\]
keeping in mind that the degrees of the variables are the negatives of
the columns of~$A$.  These true degrees are plotted as black dots in
Figure~\ref{f:primary}.  The pair of horizontal lines comes from
$\<b^2,c,d,b-e\>$, while the pair of steep diagonal lines comes from
$\<b,c^2,d,b-e\>$.
\end{example}

\begin{theorem} \label{t:sum}
Assume that $-\beta$ lies outside of the primary arrangement
$\cZ_\primary(I)$.  Then
\[
  \HH_i(E-\beta;\CC[\ttt]/I) \cong \bigoplus_{I_{\rho,J}\ \toral}
  \HH_i(E-\beta;\CC[\ttt]/\cC_{\rho,J})
\]
for all~$i$, the sum being over all toral associated primes of~$I$
from Theorem~\ref{t:primDecomp}.  In particular,
\[
  D/H_A(I,\beta) \cong \bigoplus_{I_{\rho,J}\ \toral}
  D/H_A(\cC_{\rho,J},\beta).
\]
\end{theorem}
\begin{proof}
Assume that $-\beta \notin \cZ_\primary(I)$.  Resuming the notation
from the proof of Proposition~\ref{p:Z}, we have an exact sequence $0
\to (P_I)_\toral \to P_I \to P_I/(P_I)_\toral \to 0$.  The direct sum
$\bigoplus_\andean \CC[\ttt]/\cC_{\rho,J}$ over the Andean components
of~$I$ surjects onto $P_I/(P_I)_\toral$.  Hence, by
Lemma~\ref{l:weak-rigidity}, we deduce that
$\HH_i(E-\beta;P_I/(P_I)_\toral) = 0$ for all~$i$.  Consequently,
$\HH_i(E-\beta;P_I) \cong \HH_i(E-\beta;(P_I)_\toral)$ for all~$i$.
But the latter is zero for all~$i$ by Theorem~\ref{t:rigid} because
$-\beta \notin \qdeg(P_I) \supseteq \qdeg((P_I)_\toral)$.  Therefore,
applying Euler-Koszul homology to the exact sequence in
Definition~\ref{d:P}, and using Lemma~\ref{l:weak-rigidity} to note
that this kills the Andean summands, we have proved the first display.
The second is simply the $i = 0$~case.
\end{proof}

Here is our final arrangement, outside of which the holonomic rank of
$H_A(I,\beta)$ is minimal.

\begin{definition}\label{d:jump}
Given an $A$-graded binomial~$I$, the \emph{jump arrangement}\/
of~$I$ is the union
\[
  \cZ_\jump(I) = \cZ_\andean(I) \cup \bigcup_{i=0}^{d-1}
  \qdeg\big(H^i_\mm(\CC[\ttt]/I_\toral)\big)
\]
of the Andean arrangement of~$I$ with the quasidegrees of the local
cohomology of~$\CC[\ttt]/I_\toral$ in cohomological degrees at most
$d-1$.
\end{definition}

Once the holonomic rank of a binomial $D$-module is minimal, we can
quantify it exactly.  Let $\mu_{\rho,J}$ be multiplicity
of~$I_{\rho,J}$ in~$I$ (or equivalently, in the primary component
$\cC_{\rho,J}$ of~$I$).  Denote by $\vol(A_J)$ the volume of the
convex hull of~$A_J$ and the origin, normalized so that a lattice
simplex in the group~$\ZZ A_J$ generated by the columns of~$A_J$ has
volume~$1$.

\begin{theorem}\label{t:Irank}
If $\cZ_\andean(I) \neq \CC^d$, then $H_A(I,\beta)$ has minimal rank
at~$\beta$ if and only if $-\beta$ lies outside of the jump
arrangement~$\cZ_\jump(I)$, and this minimal rank is
\[
  \rank\big(D/H_A(I,\beta)\big) = \sum_{I_{\rho,J}\text{\rm\ toral of
  dim. }d} \mu_{\rho,J} \cdot \vol(A_J).
\]
\end{theorem}
\begin{proof}
Assume that $\cZ_\andean(I) \neq \CC^d$, and denote by~$X$ the
complement of $-\cZ_\andean(I)$ in~$\CC^d$.  The global Euler-Koszul
homology $\HH_0(E-b;\CC[\ttt]/I)$ determines a sheaf on $\CC^d$, and
hence a sheaf~$\FF$ on~$X$ by restriction.  We claim that $\FF$ is a
holonomic family \cite[Definition~2.1]{MMW} over~$X$.  In fact, we
claim that $\FF$ is the restriction to~$X$ of the family determined by
$\HH_0(E-b;\CC[\ttt]/I_\toral)$, which is a holonomic family on all
of~$\CC^d$ by Theorem~\ref{t:holfamily}.  Our claim is immediate from
the sheaf (i.e., global Euler-Koszul) version
Proposition~\ref{p:Itoral}, which says that for all~$i$, if $\beta \in
X$ then $\HH_i(E-b;\CC[\ttt]/I) \cong \HH_i(E-b;\CC[\ttt]/I_\toral)$
in a neighborhood of~$\beta$.  This follows by the same proof as
Proposition~\ref{p:Itoral} itself, given the global version of
Lemma~\ref{l:weak-rigidity}.  This global version, in turn, follows
from the same proof as Lemma~\ref{l:weak-rigidity} itself with
$\beta_i$ replaced by~$b_i$ for all~$i$, the point being that $b_i =
(b_i - \beta_i) + \beta_i$ is a unit locally in~$\CC^d$ near~$\beta$,
since $b_i - \beta_i$ lies in the maximal ideal at~$\beta$.

The statement about minimality of rank is now a consequence of
Theorem~\ref{t:rankjumps} for $V = \CC[\ttt]/I_\toral$, noting that
the rank is infinite for $\beta \notin X$ by
Theorem~\ref{t:holonomic}.  To compute this minimal rank, we may
assume that $\beta$ is as generic as we like.  In particular, we
assume that $-\beta$ lies outside of the primary arrangement, and also
(by Lemma~\ref{l:dim}) outside of $\qdeg(\CC[\ttt]/\cC_{\rho,J})$ for
the components of dimension less than~$d$.  Using Theorem~\ref{t:sum},
we will be done once we show that $H_A(\cC_{\rho,J},\beta)$ has rank
$\mu_{\rho,J} \cdot \vol(A_J)$ for generic~$\beta$.

To do this, take a toral filtration of $\CC[\ttt]/\cC_{\rho,J}$.  We
are guaranteed that the number of successive quotients of dimension~$d$
is precisely the multiplicity of $I_{\rho,J}$ in~$\cC_{\rho,J}$, and
that all of the dimension~$d$ successive quotients are actually
$\ZZ^d$-translates of $\CC[\ttt]/I_{\rho,J}$ itself.  Therefore,
choosing $\beta$ to miss the quasidegree sets of the other successive
quotients, we find that the rank of $H_A(\cC_{\rho,J},\beta)$ equals
the multiplicity $\mu_{\rho,J}$ times the generic rank of
$H_A(I_{\rho,J},\beta) = H_{A_J}(I_{\rho,J},\beta)$, which is
$\vol(A_J)$ by \cite{Adolphson}.
\end{proof}

\begin{remark}\label{rk:I(B)}
If $I = I(B)$ is a lattice basis ideal (Example~\ref{ex:I(B)}), then
the sum in Theorem~\ref{t:Irank} can be simplified by gathering the
terms $\mu_{\rho,J}\cdot\vol(A_J)$ for which the domain of~$\rho$ is a
fixed toral saturated sublattice $L \subseteq \ZZ^J$.  The single term
that results is $\mu(L,J)\cdot\vol(A_J) = |L/\ZZ
B\cap\ZZ^J|\cdot\mu_{\rho,J}\cdot\vol(A_J)$, where $\rho : L \to
\CC^*$ is any partial character that is trivial on~$\ZZ B$.  Indeed,
the number of choices for $\rho$ is $|L/\ZZ B\cap\ZZ^J|$, and once
$I_{\rho,J}$ is associated to~$I(B)$, the same is true for any other
choice of~$\rho$; this is because rescaling the variables by a partial
character that is trivial on~$\ZZ B$ induces an automorphism of the
polynomial ring fixing the lattice basis ideal~$I(B)$.  For the same
reason, the multiplicities of the various choices of~$I_{\rho,J}$
in~$I(B)$ are all equal.  See Section~\ref{s:binom} for the relevance
of this simplification.
\end{remark}

\begin{remark}
The arrangement that we should require $-\beta$ to avoid for~$\beta$
to be called truly \emph{generic} is the union of the jump
arrangement~$\cZ_\jump(I)$ and the \emph{top arrangement}\/
$\cZ_{\text{\rm top}}(I) = \qdeg\big(\bigoplus_{\toral<d}
\CC[\ttt]/\cC_{\rho,J}\big)$, where the direct sum is over all toral
components of~$I$ with $\dim(\CC[\ttt]/I_{\rho,J}) \leq d-1$.  For
$-\beta \notin \cZ_\jump(I) \cup \cZ_{\text{\rm top}}(I)$, the module
$D/H_A(I,\beta)$ has minimal holonomic rank and decomposes as a direct
sum over the dimension~$d$ toral components.
\end{remark}

\begin{corollary} \label{c:homog-ranks}
If $I$ is standard $\ZZ$-graded without any Andean components, and
$\CC[\ttt]/I$ has Krull dimension~$d$, then the generic rank of
$H_A(I,\beta)$ equals the $\ZZ$-graded degree of~$I$. \qed
\end{corollary}
  
We close this section by illustrating a particular case of a Mellin
system \cite{mellin,ds06}.  Such systems arise when showing that
algebraic functions satisfy hypergeometric equations.  The goal of the
example is to give an instance when the local solution space of the binomial
$D$-module $D/H_A(I,\beta)$ for some nonzero parameter~$\beta$ fails
to split as a direct sum of the local solution spaces to 
binomial $D$-modules arising from components.  
Note that $\beta = 0$ always lies in the primary
arrangement: the residue field $\CC = \CC[\ttt]/\mm$ is a quotient of
every primary component $\CC[\ttt]/\cC_{\rho,J}$ because the
$A$-grading is positive (i.e., $\NN A$ is a pointed semigroup).

\begin{example}\label{ex:ds06}
Let
\[
  A = \left[\begin{array}{rrrr}
	1 & 1 & 1 & 1 \\
	3 & 2 & 1 & 0
      \end{array}\right] \quad\text{and}\quad
  B = \left[\begin{array}{rr}
	-2 & -1 \\ 
	 3 &  0 \\
	 0 &  3 \\ 
	-1 & -2
      \end{array}\right].
\]
In this case we have
\[
  I_{\ZZ B} = I(B) = \<\del_1^2\del_4-\del_2^3,
  \del_1\del_4^2-\del_3^3\> \subseteq
  \CC[\del_1,\del_2,\del_3,\del_4].
\]
That is, the lattice basis ideal $I(B)$ coincides with the lattice
ideal~$I_{\ZZ B}$.  The primary decomposition of $I_{\ZZ B}$ is
obtained from that of the ideal~$I$ in Examples~\ref{ex:qdeg}
and~\ref{ex:primary} by omitting the Andean component $\<a,c,d\>$ and
erasing all occurrences of $b-e$.  Thus the primary arrangement
of~$I_{\ZZ B}$ consists of the four lines in Figure~\ref{f:primary}
corresponding to toral components.

Let $\beta = -\left[\twoline {0}{1}\right]$.  The solutions of the
system $H_A\big(I_{\ZZ B},-\left[\twoline 0{1}\right]\!\big)$ are as
follows.  For $x=(x_1,x_2)$, let $z_1(t)$, $z_2(t)$ and $z_3(t)$ be
the local roots in a neighborhood of $(0,0)$~of
\[
  z^3+x_1z^2+x_2z+1=0.
\]
By \cite{Sturmfels:alg-hyp}, a local basis of solutions of the
$A$-hypergeometric system $H_A\big(\! -\left[\twoline
0{1}\right]\!\big) = \quad$ $H_A\big(I_A,-\left[\twoline 0{1}\right]\!\big)$
for the toric ideal~$I_A$~(\ref{eq:IA}) is given by the three roots of
the homoge\-ne\-ous equation
\[
  x_0 z^3+x_1 z^2 + x_2 z + x_3 = 0,
\]
and the solutions for the other two components are the roots of
\[
  x_0 z^3+ x_1 z^2+ \omega x_2 z + x_3 = 0
  \quad\text{and}\quad
  x_0 z^3+x_1z^2+\omega^2x_2z+x_3 = 0,
\]
where $\omega$ is a primitive cube root of~$1$.  The system
$H_A\big(I_{\ZZ B},-\left[\twoline 0{1}\right]\!\big)$ has nine
algebraic solutions coming from the roots $z=z(x_0, x_1,x_2,x_3)$
of the above equations.

This looks good: the quotient $\CC[\ttt]/I_{\ZZ B}$ is Cohen-Macaulay,
so $\HH_0(E-\beta,\CC[\ttt]/I_{\ZZ B})$ has holonomic rank that is
constant as a function of $\beta \in \CC^2$, by the rank minimality in
Theorem~\ref{t:Irank}, and equal to~$9$ because $\vol(A) = 3$.

However, the nine algebraic
solutions mentioned above only span a vector space
of dimension~$7$, not~$9$.  This means that there are two extra linearly
independent local solutions, which are non-algebraic;
see \cite[Example~4.2, Theorem~4.3, Example~4.4]{ds06}.

The binomial $D$-module explanation for this collapsing from
dimension~$9$ to dimension~$7$, and the concomitant extra two
logarithmic solutions, is that $-\beta = \left[\twoline 01\right] \in
\cZ_\primary(I_{\ZZ B})$; again see Figure~\ref{f:primary}.  Let us be
more precise.  The exact sequence in Definition~\ref{d:P} reads
\[
  0 \to \CC[\ttt]/I_{\ZZ B} \to R_0 \oplus R_1 \oplus R_2 \to
  P_{I_{\ZZ B}} \to 0,
\]
where $R_i = \CC[\del_1,\del_2,\del_3,\del_4]/\<\omega^i
\del_2\del_3-\del_1\del_4,\del_2^2-\omega^i\del_1\del_3,\omega^{2i}
\del_3^2-\del_2\del_4\>$.  The surjection to $P_{I_{\ZZ B}}$ factors
through the projection $R_0 \oplus R_1 \oplus R_2 \to \ol R \oplus \ol
R \oplus \ol R$, where $\ol R$ is the monomial quotient
$\CC[\ttt]/\<\del_2\del_3, \del_1\del_4, \del_2^2, \del_1\del_3,
\del_3^2, \del_2\del_4\>$, the coordinate ring of the intersection
scheme of any pair of irreducible components of the variety of~$I_{\ZZ
B}$.  The image of $\CC[\ttt]/I_{\ZZ B}$ in this projection is the
diagonal copy of~$\ol R$, so $P_{I_{\ZZ B}}$ is a direct sum $\ol R
\oplus \ol R$ of two copies of~$\ol R$.

On the other hand, each of the rings $R_i$ is also Cohen-Macaulay, so
the only nonvanishing Euler-Koszul homology of $R_0 \oplus R_1 \oplus
R_2$ is the zeroth.  Thus we have an exact sequence
\[
  0 \to \HH_1(E-\beta;P_{I_{\ZZ B}}) \to D/H_A(I_{\ZZ B},\beta) \to
  \bigoplus_{i=0}^2 \HH_0(E-\beta;R_i) \to \HH_0(E-\beta;P_{I_{\ZZ
  B}}) \to 0.
\]
In general, for $-\beta$ lying on precisely one of the four lines in
$\cZ_{\text{primary}}(I_{\ZZ B}) = \qdeg(P_{I_{\ZZ B}})$, the leftmost
and rightmost $D$-modules here have rank precisely~$2$, and this is
the~$2$ that causes the dimension collapse 
and the pair of logarithmic solutions to appear.
  
Given our choice of parameter $\beta = -\left[\twoline 0{1}\right]$,
for instance, the $\HH_1$ and the~$\HH_0$ in question are isomorphic
to one another, since both are isomorphic to a direct sum of two
copies of \mbox{$\HH_0(E - \left( - \left[\twoline 0{1}\right] \right);
\del_3\CC[\ttt]/\<\del_1,\del_2,\del_3\>)$}, where the $\del_3$ in
front of $\CC[\ttt]$ means to take an appropriate $A$-graded translate
(namely by $\deg(\del_2) = -\!\left[\twoline 11\right]$); this
corresponds to the upper of the two steep diagonal lines in
Figure~\ref{f:primary}.
\end{example}

\section{Local solutions of Horn $D$-modules}\label{s:horn}

We now return to the Horn hypergeometric $D$-modules---that is,
binomial $D$-modules arising from lattice basis ideals---that
motivated this work.  Theorem~\ref{thm:solsfromcomp}, the main result
of this section, provides a combinatorial formula for the generic rank
of a binomial Horn system by explicitly describing a basis for its
local solution space.  The basis we construct involves GKZ
hypergeometric functions.

Throughout this section, let $B$ and $A$ be integer matrices as in
Convention~\ref{conv:B}.  Since we have an explicit description for
the components of a lattice basis ideal~$I(B)$ at toral minimal
primes, namely Example~\ref{ex:I(B)'}, we make use of it to
compute---just as explicitly---the local solutions for generic~$\beta
\in \CC^d$ of the corresponding hypergeometric system.

\begin{convention}\label{conv:J}
Suppose that after permuting the rows and columns of~$B$, there
results a decomposition of~$B$ as in~(\ref{eq:MNOB}), where $M$ is a
$q \times p$ matrix of full rank~$q$.  Write $\oJ = \oJ(M)$ for the
$q$ rows occupied by~$M$ inside of~$B$ (before permuting), and let $J
= \{1,\ldots,n\} \minus \oJ$ be the rows occupied by~$B_J$.  Split the
variables $x_1,\ldots,x_n$ and $\del_1,\ldots,\del_n$ into two blocks
each:
\begin{align*}
  x_J = \{x_j : j \in J\}
  &\quad\text{and}\quad
  x_\oJ = \{x_j : j \notin J\}.
\\
  \ttt_J = \{\del_j : j \in J\}
  &\quad\text{and}\quad
  \ttt_\oJ = \{\del_j : j \notin J\}.
\end{align*}
As before, $A_J$ is the submatrix of~$A$ with columns $\{a_j : j \in J\}$.
\end{convention}

With the notation above, fix for the remainder of this article a toral
prime $I_{\rho,J}$ of $I(B)$. Since $I(B)$ is generated by $m=n-d$
elements, $I_{\rho, J}$ has dimension at least $d$.  On the other
hand, toral primes can have dimension at most $d$, by
Lemma~\ref{l:dim}.  Thus we have the following.

\begin{lemma} \label{lemma:d+m}
All toral primes of the lattice basis ideal $I(B)$ have dimension
exactly $d$ and are minimal primes of $I(B)$.\qed
\end{lemma}

\begin{observation}\label{obs:square}
Since the dimension of $I_{\rho,J}$ equals $n-p-(m-q) = d+q-p$, if
$I_{\rho,J}$ is toral, then the previous lemma implies that $q=p$.
Thus, from now on, the matrix $M$ is a $q\times q$ mixed invertible
matrix (and $q$~is allowed to be~$0$).
\end{observation}

Recall from Example~\ref{ex:I(B)'} that for a toral minimal
prime~$I_{\rho,J}$, the component can be written as $\cC_{\rho,J} =
I(B) + I_\rho + U_M$, where $U_M \subseteq \CC[\ttt_\oJ]$ is the
ideal $\CC$-linearly spanned by all monomials with exponent vectors in
the union of the unbounded $M$-subgraphs of~$\NN^\oJ$
(Definition~\ref{d:M}).

In order to construct local solutions of $H_A(\cC_{\rho,J},\beta)$ we
need two ingredients: local solutions of a GKZ-type system
$H_{A_J}(I_{\rho},\beta')$ and polynomial solutions of the constant
coefficient system~$I(M) = \<\ttt^u - \ttt^v : u - v \text{ is a
column of } M, \, u, v \in \NN^n\>$.  As it turns out, solving the
differential equations $I(M)$ is equivalent to finding the
$M$-sub\-graphs of $\NN^\oJ$.

\begin{lemma}\label{lemma:pre-crucial}
Let $M$ be a $q\times q$ mixed invertible integer matrix, and assume
that $q>0$.  Fix $\gamma \in \NN^\oJ$, and denote by~$\Gamma$ the
$M$-subgraph containing~$\gamma$.
\begin{numbered}
\item
The system $I(M)$ of differential equations has a unique formal power
series solution of the form $G_{\gamma} = \sum_{u \in \Gamma}
\lambda_u x^u$ in which $\lambda_{\gamma} = 1$.
\item
The other coefficients $\lambda_u$ of $G_{\gamma}$ for $u \in
\Gamma$ are all nonzero.
\end{numbered}
\end{lemma}

This lemma will be proved together with Proposition~\ref{prop:crucial}.

\begin{notation}\label{not:representatives}
Given a $q \times q$ mixed invertible matrix $M$, we fix a 
set $\cS(M) \subset \NN^{\oJ}$ of representatives for the bounded 
$M$-sugraphs of $\NN^\oJ$. 
In particular, the cardinality of $\cS(M)$ equals the number of
bounded $M$-subgraphs, which we denote by $\mu_M$.  If $q=0$, we set
$\cS(M) = \{ \nothing \}$ and declare $\mu_M$ to be~$1$.
\end{notation}

\begin{prop}\label{prop:crucial}
With the notation from Lemma~\ref{lemma:pre-crucial} and
Notation~\ref{not:representatives},
\begin{numbered}
\item
The set $\{G_{\gamma} : \gamma$ runs over a set of representatives for
the $M$-subgraphs of\/~$\NN^\oJ\}$ is a basis for the space of all
formal power series solutions of~$I(M)$.
\item
The set $\{G_{\gamma} : \gamma \in \cS(M) \}$ is a basis for the space
of polynomial solutions of~$I(M)$.
\end{numbered}
\end{prop}

\begin{proof}[Proof of Lemma~\ref{lemma:pre-crucial} and
Proposition~\ref{prop:crucial}]

We begin with the first statement from
Lem\-ma~\ref{lemma:pre-crucial}.  If\/ $\Gamma = \{ \gamma \}$ then
$G_{\gamma} = x^{\gamma}$.  We check that this is a solution of $I(M)$
working by contradiction.  Let $w$ be a column of $M$ such that
$\ttt^{w_+} x^{\gamma} \neq \ttt^{w_-} x^{\gamma} $.  Then one
of these terms is nonzero, say $\ttt^{w_+} x^{\gamma}$, so that
$\gamma-w_+ \in \NN^\oJ$.  But then $\gamma-w_+ + w_- = \gamma - w \in
\NN^\oJ$, and so $\gamma - w \in \Gamma$, a contradiction, because
$\gamma - w \neq \gamma$ and $\Gamma$ is a singleton.

Now assume that $\Gamma$ is not a singleton, and fix $u \in \Gamma$
such that $u - \gamma = w$ is a column of~$M$.  We want to define the
coefficients of $G_{\gamma}$, and we will start with~$\lambda_u$.
Since $u - \gamma = w = w_+ - w_-$, we have $u - w_+ = \gamma - w_-
\in \NN^\oJ$, since $u$ and $\gamma$ both lie in~$\NN^\oJ$ and the
supports of $w_+$ and~$w_-$ are disjoint.  Set $\lambda_u =
\ttt^{w_-}(x^\gamma)/\ttt^{w_+}(x^u)$, and observe that numerator and
denominator are \emph{nonzero} constant multiples of $x^{u-w_+} =
x^{\gamma-w_-}$.  Use this procedure to define the coefficients
corresponding to the neighbors of $\gamma$.  Now, if we
know~$\lambda_u$ and we are given a neighbor $u'\in \Gamma$ of~$u$,
say $u' - u = w$, then set $\lambda_{u'} = \ttt^{w_-}(x^u)/\lambda_u
\ttt^{w_+}(x^{u'})$.  Propagating this procedure along~$\Gamma$ we
obtain all of the coefficients~$\lambda_u$.  The formal power series
$G_{\gamma}$ defined this way is tailor-made to be a solution
of~$I(M)$.

Since $M$-subgraphs are disjoint, it is clear that the series
$G_{\gamma}$ are linearly independent. Now let $G = \sum_{u \in
\NN^\oJ} \nu_u x^u$ be a formal power series solution of $I(M)$.  We
claim that $G-\nu_{\gamma}G_{\gamma}$ has coefficient zero on all
monomials from~$\Gamma$.  This follows from the fact that
$G-\nu_{\gamma}G_{\gamma}$ has coefficient zero on the
monomial~$x^\gamma$; indeed, if the difference contained a monomial
from~$\Gamma$, it would have to contain~$x^\gamma$ with a nonzero
coefficient, as can be seen by the propagation argument from before.
(The uniqueness of~$G_{\gamma}$ that we need for
Lemma~\ref{lemma:pre-crucial} also follows from this argument.)  It is
now clear that our candidate power series solution basis is a spanning
set, and the statement for polynomial solutions has the same proof.
\end{proof}

\begin{remark}
The system $I(M)$ is itself a binomial Horn system; there are no Euler
operators because $M$ is invertible.  We stress that it is a very
special feature of hypergeometric differential equations that their
irreducible (Puiseux) series solutions are determined (up to a
constant multiple) by their supports.  In general, this is far from
being the case for systems of differential equations that are not
hypergeometric.
\end{remark}

We can use this correspondence between $M$-subgraphs and solutions of
$I(M)$ to compute examples.

\begin{example}
\label{ex:concrete-subgraph}
Consider the $3\times 3$ matrix
\[ M = \left[ \begin{array}{rrr}
 1 & -5 &  0 \\
-1 &  1 & -1 \\
 0 &  3 &  1
\end{array}
\right]
\]
A basis of solutions (with irreducible supports) of
$I(M)$ is easily computed:
\[
\left\{
  1, \quad x+y+z, \quad  (x+y+z)^2, \quad (x+y+z)^3, \quad
  \sum_{n \geq 4} \frac{(x+y+z)^n}{n!}
\right\}.
\]
\begin{figure} 
\[
\includegraphics{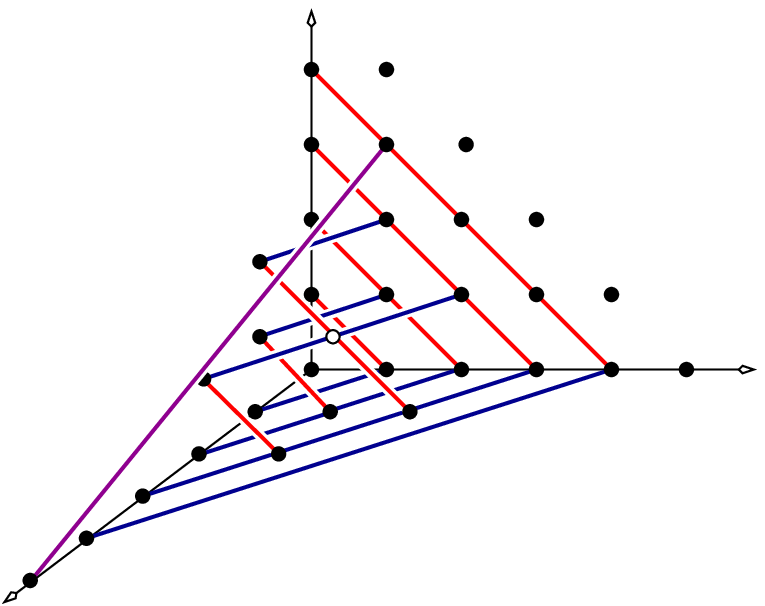}
\]
\caption{The $M$-subgraphs of $\NN^3$}
\label{f:M-subgraphs}
\end{figure}

The $M$-subgraphs of $\NN^3$ are the four slices $\{(a,b,c) \in \NN^3:
a+b+c=n\}$ for $n \leq 3$; for $n \geq 4$, two consecutive slices are
$M$-connected by $(-5,1,3)$, yielding one unbounded $M$-subgraph.
\end{example}

The following definition will allow us to determine a set of
parameters $\beta$ for which the system $H_A(\cC_{\rho,J},\beta)$ has
the explicit basis of solutions that we construct for
Theorem~\ref{thm:solsfromcomp}.

\begin{definition}
A \emph{facet} of $A_J$ is a subset of its columns that is maximal
among those minimizing nonzero linear functionals on~$\ZZ^d$.  For a
facet $\sigma$ of~$A_J$ let $\nu_\sigma$ be its \emph{primitive
support function}, the unique rational linear form satisfying
\begin{enumerate}
\item $\nu_\sigma(\ZZ A_J) = \ZZ$,
\item $\nu_\sigma(a_j) \geq 0$ for all $j \in J$,
\item $\nu_\sigma(a_j) = 0$ for all $a_j \in \sigma$.
\end{enumerate}
A parameter vector $\beta \in \CC^d$ is \emph{$A_J$-nonresonant}\/ if
$\nu_\sigma(\beta) \notin \ZZ$ for all facets $\sigma$ of~$A_J$.  Note
that if $\beta$ is $A_J$-nonresonant, then so is $\beta+A_J(\gamma)$
for any $\gamma \in \ZZ^J$.
\end{definition}

The reason nonresonant parameters are convenient to work with is the
following.

\begin{lemma}\label{l:iso}
If $\beta$ is $A_J$-nonresonant, then for any $\gamma \in \NN^J$, and
for all torus translates $I_\rho$ of the toric ideal~$I_{A_J}$, right
multiplication by $\ttt_J^\gamma$ induces a $D_J$-module
isomorphism $D_J/H_{A_J}(I_\rho,\beta) \to
D_J/H_{A_J}(I_\rho,\beta+A_J(\gamma))$, whose left inverse we denote
by $\ttt_J^{-\gamma}$.
\end{lemma}
\begin{proof}
For $R = \CC[\ttt_J]/I_\rho$ there is an exact sequence $0 \to R
\overset{{\ttt_J^{\gamma}}}{\longrightarrow} R \to R/\ttt_J^\gamma R
\to 0$.  Since the multiplication by $\ttt_J^{\gamma}$ occurs in the
right-hand factor of $D_J \otimes_{\CC[\ttt_J]} R$, the map on
Euler-Koszul homology over~$D_J$ induced by $\ttt_J^\gamma$
corresponds to right multiplication.  But $R/\ttt_J^{\gamma}R$ is
toral by Lemma~\ref{l:W}, and its set of quasidegrees is the Zariski
closure of $\{ - A_J \vartheta : \vartheta \in \NN^{J}, \vartheta_i <
\gamma_i \; \mbox{for some} \; i \in J\}$, which is a finite subspace
arrangement contained in the resonant parameters.  Now apply
Lemma~\ref{l:V} and Theorem~\ref{t:rigid} to complete the proof.
\end{proof}

The following definition characterizes parameter vectors with
particularly nice behavior when it comes to isomorphisms between
$H_A(I,\beta)$ for varying $\beta$.

\begin{definition}\label{def:verygen}
A parameter vector $\beta \in \CC^d$ is called \emph{very generic} if
$\beta - A_\oJ(\gamma)$ is $A_J$-nonresonant for every \mbox{$\gamma
\in \cS(M)$}.
\end{definition}

\begin{remark}\label{rem:integer-derivatives}
Denote by ${\rm Sol}(I_{\rho},\beta)$ the space of local holomorphic
solutions of $H_{A_J}(I_{\rho},\beta)$ near a nonsingular point.
Given $\alpha \in \NN^J$, the $D$-module isomorphism in
Lemma~\ref{l:iso} induces a vector space isomorphism
\[
  {\rm Sol}(I_{\rho},\beta)
  \longleftarrow
  {\rm Sol}(I_{\rho},\beta+A_J\alpha)
\]
given by differentiation by $\ttt_{J}^{\alpha}$. If we denote the
inverse of this map by $\ttt_{J}^{-\alpha}$, a number of questions
arise: for instance, given a local solution $f \in {\rm
Sol}(I_{\rho},\beta)$ where $\beta$ is very generic, 
and taking for instance $J = \{1,2\}$,
\begin{itemize}
\item
is $\ttt_{\{1,2\}}^{-(1,0)} \big(\ttt_{\{ 1,2\}}^{-(0,1)} f \big)$ 
equal to $\ttt_{\{1,2\}}^{-(1,1)}f$?
\item
is $\ttt_{\{1,2\}}^{(1,1)} \big(\ttt_{\{ 1,2\}}^{-(2,2)} f\big)$ equal
to $\ttt_{\{1,2\}}^{-(1,1)} f$?
\end{itemize}
Both questions have positive answers; their verification is based on
the fact that the left and right inverses of a vector space
isomorphism are the same. We conclude that $\ttt_J^{-\alpha}f$ is
well-defined for any $f \in {\rm Sol}(I_{\rho},\beta+A_J\alpha)$, if
$\beta$ is very generic and $\alpha$ is an arbitrary integer~vector.

At the level of $D$-modules, however, $\ttt_J^{-\alpha}$ for $\gamma
\in \ZZ^{J}$ is not necessarily well-defined, because the right and
the left inverses of a $D$-isomorphism need not coincide.
\end{remark}

We resume Notation~\ref{not:representatives}.  If $q=0$, then set
$G_{\nothing}=1$.  If $q>0$ and $\gamma \in \cS(M)$, then rewrite the
polynomial $G_\gamma$ from Lemma~\ref{lemma:pre-crucial} as follows:
\[
  G_{\gamma} = x_\oJ^\gamma \sum_{\gamma+Mv\in \Gamma} c_v x_\oJ^{Mv}.
\]
By Proposition~\ref{prop:crucial}, $\{G_{\gamma} : \gamma \in
\cS(M)\}$ is a basis for the polynomial solution space of~$I(M)$.

Given a local solution $f = f(x_J)$ of the system
$H_{A_J}(I_{\rho},\beta -A_\oJ(\gamma))$ for some $\gamma \in \cS(M)$,
define
\begin{equation}\label{eqn:solform}
  F_{\gamma,f} = x_\oJ^\gamma \sum_{\gamma+Mv \in \Gamma} c_v x_\oJ^{Mv}
  \ttt_J^{-Nv}(f),
\end{equation}
where $\ttt_J^{-Nv}f$ is as in Remark~\ref{rem:integer-derivatives}.
Note that if $q=0$, we have $F_{\nothing,f}=f$.

The condition of being very generic is open and dense in 
the standard topology of~$\CC^d$, so
that the rank of $H_A(\cC_{\rho,J},\beta)$ for such parameters equals
the generic rank of this binomial $D$-module, in the sense of
Theorem~\ref{t:Irank}.

\begin{theorem}\label{thm:solsfromcomp}
Let $\cC_{\rho,J}$ be a toral component of~$I(B)$ and let $\beta \in
\CC^d$ be a very generic parameter vector.  Given $\gamma \in \cS(M)$,
fix a basis $\cB_{\gamma}$ of local solutions of $H_{A_J}(I_\rho,\beta
- A_\oJ(\gamma))$.  The $\mu_M \cdot \vol(A_J)$ functions $\{
F_{\gamma,f} : \gamma \in \cS(M), f \in \cB_{\gamma}\}$ form a local
basis for the solution space of the binomial
$D$-module~$D/H_A(\cC_{\rho,J},\beta)$.
\end{theorem}

Before proving Theorem~\ref{thm:solsfromcomp}, let us see the
construction~(\ref{eqn:solform}) in some explicit examples.

\begin{example}
Consider the matrices
\[
  A =
  \left[\begin{array}{rrrrr}
	1 &  1 & 1 & 1 & 1 \\
	5 & 10 & 0 & 7 & 6 
  \end{array}\right]
\quad\text{and}\quad
  B =
  \left[\begin{array}{rrr}
	 0 & -1 &  2 \\
	-1 &  0 & -1 \\
	 0 &  1 & -1 \\
	 4 &  5 &  0 \\
	-3 & -5 &  0
  \end{array}\right].
\]
We concentrate on the decomposition
\[
  M =
  \left[\begin{array}{rr}
	 4 &  5 \\
	-3 & -5
  \end{array}\right];
\quad
  N =
  \left[\begin{array}{rr}
	 0 & -1 \\
	-1 &  0 \\
	 0 &  1
  \end{array}\right]; 
\quad
 B_J =
  \left[\begin{array}{r}
	 2 \\
	-1 \\
	-1
  \end{array}\right].
\]
Note that $\ZZ B_J$ is saturated, so there is only one associated
prime coming from this decomposition, namely $I_{\left[
\begin{smallmatrix} 1 & 1 & 1 \\ 5 & 10 & 0 \end{smallmatrix}\right]}
+ \<\del_4,\del_5\>$, and this is toral since $\det(M)\neq 0$.

The polynomial $\varphi = 5 x_4^4 x_5^2+ 2 x_4^5+ 2 x_5^5+ 40x_4x_5^3$
is a solution of the constant coefficient system $I(M)$.  Let $f$ be a
local solution of the $\left[\begin{smallmatrix} 1 & 1 & 1 \\ 5 & 10 & 0
\end{smallmatrix}\right]$-hypergeometric system that is homogeneous
of degree $\beta - \left[ \begin{smallmatrix} 6 \\ 40
\end{smallmatrix}\right]$.  It can be verified that the following
function is a solution of $H(B,\beta)$:
\[
  5 x_4^4 x_5^2 f + 
  2 x_4^5 \del_1^{1} \del_2^{-1} \del_3^{-1} f + 
  2 x_5^5 \del_2^{-1} f+ 
  40 x_4x_5^3 \del_1 \del_2^{-2}\del_3^{-1} f.
\]
In this example, the new solution we constructed has $1$-dimensional
support.
\end{example}

\begin{example}
Our procedure for constructing solutions works even when $M$ is an $m
\times m$ matrix, i.e., $M$ is a maximal square submatrix of $B$.  For
instance, consider
\[
  A =
  \left[\begin{array}{rrrrr}
	1 & 1 & 1 & 1 \\
	1 & 0 & 3 & 2
  \end{array}\right]
\quad\text{and}\quad
  B =
  \left[\begin{array}{rr}
	 2 & -3 \\
	-1 &  2 \\
	 0 &  1 \\
	-1 &  0
  \end{array}\right].
\]
We concentrate on the component
\[
  M =
  \left[\begin{array}{rr}
	 2 & -3 \\
	-1 &  2
  \end{array}\right];
\quad
  N =
  \left[\begin{array}{rr}
	 0 & 1 \\
	-1 & 0
  \end{array}\right];
\quad
  B_J =
  \nothing.
\]
Again, we only have one (toral) component, associated to
$\<\del_1,\del_2\>$.  Let $p = x_1^2 + 2x_2$.  This is a solution of
$I(M) = \<\del_1^2-\del_2, \del_1^3-\del_2^2\>$.  Since $B_J$ is
empty, we need only consider solutions of the homogeneity equations
that are functions of $x_3$ and~$x_4$.  Since $\det(M) = 1 \neq 0$,
the complementary minor of~$A$ is also nonzero, and therefore there
exists a unique monic monomial in $x_3$ and~$x_4$ of each degree.  To
make a solution of $I + \<E - \beta\>$, let $x_3^{w_3}x_4^{w_4}$ be
the unique monic solution of the homogeneity equations with parameter
$\beta - \left[\twoline 10\right]$.  Then
\[
  w_4 x_1^2 x_3^{w_3} x_4^{w_4-1}  + 2 x_2  x_3^{w_3}x_4^{w_4} =
  x_3^{w_3}x_4^{w_4-1}(w_4x_1^2+2x_2)
\]
is the desired solution of $H(B,\beta)$.
\end{example}

\begin{proof}[Proof of Theorem~\ref{thm:solsfromcomp}]
First note that if $q=0$, all of the statements hold by construction.
Therefore we assume that $q \geq 2$.

It is clear that all of the binomial generators of $I_\rho$
annihilate~$F_{\gamma,f}$.  It is also easy to check that
$F_{\gamma,f}$ satisfies the desired homogeneity equations.  Let then
$(\nu, \delta) \in \ZZ^J \times \ZZ^\oJ$ be one of the $q$ columns
of~$B$ involving $N$ and~$M$; i.e., $(\nu, \delta) = \left[\twoline
NM\right]\!e_k$ for some $k \in \oJ$.  To prove that
$(\ttt_J^{\nu_+}\ttt_\oJ^{\delta_+} -
\ttt_J^{\nu_-}\ttt_\oJ^{\delta_-})(F_{\gamma,f}) = 0$, notice
that $(\ttt_\oJ^{\delta_+}-\ttt_\oJ^{\delta_-})(G_{\gamma}) =
0$, which implies that for all $v$ with $c_v \neq 0$, either
$\ttt_\oJ^{\delta_+}(x^{\gamma+Mv}) = 0$
or there exists another integer vector $w$ with $c_w \neq 0$ such that
\[
  \ttt_\oJ^{\delta_+} (c_v x_\oJ^{\gamma + Mv)}) =
  \ttt_\oJ^{\delta_-} (c_w x_\oJ^{\gamma + Mw}).
\]
In the first case, $\ttt_J^{\nu_+}\ttt_\oJ^{\delta_+} \left(
c_v x_\oJ^{\gamma + Mv} \ttt_J^{N(-v)}(f) \right) = 0$.  In the
second case, since the monomials on the left and right-hand sides of
the above equation must have the same exponent vector, we see that
$\delta = M(v-w) = M \cdot e_k$.  But $M$ is invertible by assumption,
so that $v-w = e_k$.  This implies that $\nu = Ne_k = N(v-w)$.

Consequently, $\ttt_J^{N(-v) + \nu_+}(f) = \ttt_J^{N(-w) +
\nu_-}(f)$, 
and thus
\[
  \ttt_\oJ^{\delta_+}\ttt_J^{\nu_+} \left(c_v x_\oJ^{\gamma + Mv}
  \ttt_J^{N(-v)}(f)\right)
  = \ttt_\oJ^{\delta_-}\ttt_J^{\nu_-} \left(c_w x_\oJ^{\gamma + Mw}
  \ttt_J^{N(-w)}(f) \right).
\]
Moreover, it is clear that the $F_{\gamma,f}$ are linearly independent.

Now we need to show that these functions span the local solution space
of $H_A(\cC_{\rho,J},\beta)$.  Let $F=F(x_1,\dots,x_n)$ be a local
solution of $H_A(\cC_{\rho,J},\beta)$.  Here we use the explicit
description of $\cC_{\rho,J}$ from Example~\ref{ex:I(B)'}.  Since the
monomials in~$U_M$ annihilate~$F$, we can write
\[
  F=\sum_{\gamma \in \cS(M)} \sum_{\alpha \in \Gamma} x_\oJ^\alpha
  h_\alpha(x_J),
\]
where the sum runs over 
the union of the bounded $M$-subgraphs, that is, the sum runs over all
$\alpha$ such that $\ttt_J^\alpha$ does not belong to~$U_M$.

The functions $h_\alpha$ are solutions of $H_A(I_\rho,\beta - A_\oJ
(\alpha))$, as is easy to check.  Note that $h_\alpha$ may be zero.

Now it is time to use the equations~$I(B)$.  First, observe that we
may assume that the $x_\oJ$-monomials in $F$ belong to a single
$M$-subgraph of $\NN^\oJ$.  This is because the only equations
relating different summands from~$F$ are those from $I(B)$, which will
relate a summand $x_\oJ^\alpha h_\alpha(x_J)$ to a different summand
$x_\oJ^{\alpha'} h_{\alpha'}(x_J)$ exactly when $\alpha-\alpha'$ or
$\alpha'-\alpha$ is a column of~$M$, thus staying within an
$M$-subgraph.

So fix a bounded $M$-subgraph $\Gamma$ corresponding to a $\gamma \in
\cS(M)$, and write
\[
  F=\sum_{\alpha \in \Gamma} x_\oJ^\alpha h_\alpha(x_J).
\]
Fix $\alpha \in \Gamma$ such that $h_\alpha \neq 0$, recall that
$G_{\gamma}$ is a polynomial solution of $I(M)$ whose support is
$\Gamma$, and let $c$ be the (nonzero) coefficient of $x_\oJ^\alpha$
in~$G_{\gamma}$.  We want to show that $F = (1/c) F_{\gamma,
h_\alpha}$.  Since we know that $F-(1/c)F_{\gamma, h_\alpha}$ has
support contained in~$\Gamma$ and has no summand with $x^\alpha$, the
desired equality will be a consequence of the following.

\textbf{Claim.} With the notation above, if $h_\alpha=0$ then $F = 0$.

{\it Proof of the Claim.}
If $\Gamma$ is a singleton, we are done. Otherwise pick $\alpha' \in
\Gamma$ such that $\alpha'-\alpha$ or $\alpha -\alpha'$ is a column of
$M$, say $\alpha - \alpha' = M e_k$.  The binomial from the
corresponding column of $B$ is $\ttt_J^{{Ne_k}_+}
\ttt_\oJ^{{Me_k}_+} - \ttt_J^{{Ne_k}_-} \ttt_\oJ^{{M
e_k}_-}$.  Since this binomial annihilates $F$, and $\alpha-(M e_k)_+
= \alpha' - (M e_k)_-$, we have
\[
  \ttt_\oJ^{(M e_k)_+}x_\oJ^\alpha \ttt_J^{(N e_k)+}h_\alpha 
  =
  \ttt_J^{(M e_k)_-}x_\oJ^{\alpha'} \ttt^{(N e_k)_-}h_{\alpha'},
\]
so that, as $h_\alpha = 0$,
\[
  \ttt_\oJ^{(M e_k)_-}x_\oJ^{\alpha'} \ttt_J^{(N e_k)_-}h_{\alpha'} = 0.
\]
Now, the first derivative in the previous expression is nonzero, so
$\ttt_J^{(N e_k)_-}h_{\alpha'} = 0 $.  But then $h_{\alpha'} = 0$,
since differentiation in any of the $x_J$ variables is an isomorphism
(which is why we need our parameter to be very generic).

Propagate the previous argument along $\Gamma$ to finish the proof of
the claim, and with it the proof of the theorem.
\end{proof}

\begin{remark}
When $\cC_{\rho,J}$ is Andean (and $\beta$ is generic), the above
procedure produces no nonzero solutions, as expected, since in this
case, $D/H_A(\cC_{\rho,J},\beta) = 0$ for generic $\beta$.  The reason
that the construction breaks in this situation is that there are no
nonzero solutions for the ``toric'' part.
\end{remark}

\begin{corollary}\label{cor:formula}
Fix $B$ as in Convention~\ref{conv:B}.  If there exists a
parameter~$\beta$ for which the binomial Horn system $H(B,\beta)$ has
finite rank, then for generic parameters~$\beta$, this rank is
\[
  \rank(H(B,\beta))
  =
  \sum_{\cC_{\rho,J}\ \toral}
  \mu_M \cdot \vol(A_{J})
  =
  \sum_{B\,=\,\left(\begin{smallmatrix}
                    \!N&\!B_J\\
                    \!M&\!0
                    \end{smallmatrix}\right)}
  \mu_M \cdot g(B_J) \cdot \vol(A_J),
\]
the former sum being over all toral components $\cC_{\rho,J}$ of the
lattice basis ideal~$I(B)$, and the latter sum being over all
decompositions of~$B$ as in~(\ref{eq:MNOB}) with $M$
invertible.
Here, $g(B_J)$ is the
cardinality of\/ $\sat(\ZZ B_J)/\ZZ B_J$, and $\mu_M$ is the number of
bounded $M$-subgraphs of\/ $\NN^\oJ$.
\end{corollary}
\begin{proof}
The first equality is a direct consequence of
Theorem~\ref{thm:solsfromcomp} and Theorem~\ref{t:sum}.  Comparing
with Theorem~\ref{t:Irank} yields the fact that $\mu_M$ equals the
multiplicity $\mu_{\rho,J}$ of $I_{\rho,J}$ as an associated prime of
$I(B)$.  For the second equality, the number of components arising
from a decomposition~(\ref{eq:MNOB}) is~$g(B_J)$
\cite[Corollary~2.5]{binomialideals}.
\end{proof}

\begin{remark}
The only sense in which our rank formula for Horn systems is not
completely explicit is that it lacks an expression for the number
$\mu_M$ of bounded $M$-subgraphs.  In the case that $I(M)$ (or $I(B)$)
is a complete intersection, Cattani and Dickenstein
\cite{cattani-dick:binomial-ci} can be applied to provide an explicit
recursive formula for $\mu_M=\mu_{\rho,J}$.  The general case---even
just the toral case---of this computation is an open problem.
\end{remark}

\begin{example}
\label{ex:himalayan}
The existence clause for~$\beta$ in Corollary~\ref{cor:formula} is
essential: there exist matrices~$B$ for which holonomicity of the Horn
system $H(B,\beta)$ fails for all parameters~$\beta$.  Let
\[
  A = \left[
      \begin{array}{rrrrr}
      -3 & -1 & 2 & 1 & 0 \\
      -1 & 0 & 1 & 1 & 1
      \end{array}
      \right]
  \quad\text{and}\quad
  B = \left[
      \begin{array}{rrr}
       1 & 1 & 1 \\
      -1 & -2 & -3 \\
       1 & 0 & 0 \\
       0 & 1 & 0 \\
       0 & 0 & 1
      \end{array}
      \right].
\]
Then $\<\del_1,\del_2\>$ is an Andean prime of $I(B)$.  The
quasidegree set of the corresponding component is $\CC \cdot \left[
\begin{smallmatrix} 2 & 1 & 0 \\ 1 & 1 & 1
\end{smallmatrix}\right] = \CC^2$, which means that 
the Andean arrangement of~$I(B)$ equals~$\CC^2$, and thus $H(B,\beta)$
is non-holonomic for all parameters~$\beta$.
\end{example}

A sufficient condition to guarantee holonomicity of $H(B,\beta)$ for
generic parameters is to require that $I(B)$ be a complete
intersection.  This is automatic for $m=2$, so the following result is
a direct generalization of~\cite[Theorem 8.1]{dms}.

\begin{prop}
If $I(B)$ is a complete intersection, then the binomial Horn system
$H(B,\beta)$ is holonomic for generic parameters $\beta$.
\end{prop}
\begin{proof}
If $I(B)$ is a complete intersection, its associated primes all have
dimension $d$. In combinatorial terms, we encounter only square
matrices $M$ in the primary decomposition of $I(B)$.  The component
associated to a decomposition (\ref{eq:MNOB}) is Andean exactly when
$\det(M) = 0$, and in this case, the corresponding quasidegree set is
$\CC A_J \subsetneq \CC^d$, as $A_J$ does not have full rank. We
conclude that the Andean arrangement of $I(B)$ is strictly
contained~in~$\CC^d$.
\end{proof}

When $I(B)$ is standard $\ZZ$-graded and has no Andean components, 
we can obtain a cleaner rank formula, by noting that the sum in
Corollary~{\ref{cor:formula}} equals the degree of~$I(B)$.
This is a generalization of a result in \cite{timur}.

\begin{corollary}\label{c:noAndean}
Assume that $I(B)$ is standard $\ZZ$-graded and has no Andean
components.  Let $d_1,\dots,d_m$ be the degrees of the generators
of~$I(B)$.  Then
\[ 
  \rank(H(B,\beta)) = d_1\cdots d_m \quad \text{for all }\beta \in
  \CC^d.
\]
\end{corollary}
\begin{proof}
Since $I(B)$ has no Andean components, $\CC[\ttt]/I(B)$ is toral.
Moreover, $\CC[\ttt]/I(B)$ is Cohen-Macaulay, as $I(B)$ is a complete
intersection by Lemma~\ref{lemma:d+m}.  Theorem~\ref{t:rankjumps}
implies that the holonomic rank of $H(B,\beta)$ is constant.  Now
apply Corollary~\ref{c:homog-ranks}.
\end{proof}

In the standard $\ZZ$-graded case, binomial $D$-modules are regular
holonomic.  The method of canonical series solutions \cite{SST} then
produces expansions for their solutions into power series with
logarithms.  This method applies to any regular holonomic $D$-ideal,
not just those of the form $H_A(I,\beta)$ for $\ZZ$-graded~$I$.
However, in the binomial $D$-module case, for very generic parameters
the \emph{supports} of the series solutions (i.e., the sets of
exponents of the monomials appearing with nonzero coefficients in the
series) can be very explicitly described, owing to the fact that such
combinatorial descriptions exist for GKZ functions (see \cite{GKZ} or
\cite{SST}).

\begin{definition}\label{def:fully-supported}
Let $L \subseteq \ZZ^n$ be a rank~$m$ lattice and $\alpha \in \CC^n$.
A formal series $x^\alpha \sum_{u \in L} c_u x^u$ is \emph{fully
supported}\/ on $\alpha + L$ if there exists an $m$-dimensional
polyhedral cone $C \subseteq \RR^n$, a vector $\lambda \in L$, and a
sublattice $L' \subseteq L$ of full rank~$m$ such that every term $c_u
x^{\alpha + u}$ for $u \in (\lambda + C) \cap L'$ has nonzero
coefficient $c_u$.
\end{definition}

In our case, the lattice $L$ comes from a toral component
$\cC_{\rho,J}$ of a $\ZZ$-graded lattice basis ideal~$I(B)$
corresponding to a decomposition~(\ref{eq:MNOB}), and sublattices $L'$
are necessary because $L$ is often the saturation of some other given
lattice (such as $\ZZ B$).  Recall that $\rho : L \rightarrow \CC^*$
is a partial character, where the lattice $L \subseteq \ZZ^J \subseteq
\ZZ^n$ is the saturation of the integer span of the columns of~$B_J$
in the notation from~(\ref{eq:MNOB}); equivalently, $L =
\ker_\ZZ(A_J)$ by Lemma~\ref{l:toralAndean}.

\begin{corollary}\label{cor:supports}
Fix $\gamma\in \cS(M)$ as in Notation~\ref{not:representatives}, and
let $\Lambda \subseteq \NN^n$ be a $B$-subgraph whose projection to
the $\oJ$ coordinates is the bounded $M$-subgraph containing~$\gamma$.
If $\beta \in \CC^d$ is very generic, then there exists a vector
$\alpha \in \CC^J \subseteq \CC^n$, unique modulo~$L$, such that
$A(\alpha+\Lambda)=\beta$.  The $\vol(A_J)$ functions $\{ F_{\gamma,f}
: f \in \cB_{\gamma}\}$ from Theorem~\ref{thm:solsfromcomp} are fully
supported on $\alpha + \Lambda+ L$.
\end{corollary}
\begin{proof}
First note that $A\Lambda$ is a well defined point in $\ZZ^d$, as two
elements of $\Lambda$ differ by an element of $\ZZ B \subseteq
\ker_{\ZZ}(A)$.  It follows that $\Lambda + L = \lambda + L$ for any
$\lambda \in \Lambda$, so it makes sense to be fully supported on
$\alpha + \Lambda+ L$.  On the other hand, the linear system $A \alpha
= - A \Lambda + \beta$ has a unique solution modulo (the complex span
of)~$L$ since $A_J$ has full rank; recall that we are working with a
toral component.  Now the statement about the supports follows from
Theorem~\ref{thm:solsfromcomp}, since elements of $\cB_\gamma$ are
expressible as series on $\alpha + \Lambda + L$ that are fully
supported---either as Gamma series \`a la \cite{GKZ} or as canonical
series \`a la \cite{SST}.
\end{proof}

\begin{remark}
We saw in the Introduction that a solution of $H_A(\beta)$ (or any of
the binomial $D$-modules arising from a torus translate of~$I_A$) is
essentially a function in $m=n-d$ variables.  In fact, for
generic~$\beta$, if we choose canonical series expansions as
in~\cite{SST}, then their supports are translates cones of
dimension~$m$.  This implies that the support of the
series~(\ref{eqn:solform}) has dimension $m-q$, since the dimension of
the support equals that of any series expansion of~$f$.  In fact, this
support might not be the set of lattice points in a cone, but in a
polyhedron whose recession cone has the correct dimension.
Nonetheless, the only fully supported solutions of
$H_A(I(B),\beta)=H(B,\beta)$ arise from $H_A(I_{\ZZ B},\beta)$.
Interestingly, there can be no solutions with support of
dimension~\mbox{$m-1$}, because a matrix with $q = 1$ row is never
mixed.  This explains why Erd\'elyi only found Puiseux polynomial
solutions (such as in Examples~\ref{ex:erdelyi}, \ref{ex:erdelyi'},
\ref{ex:erdelyi''}, and~\ref{ex:erdelyi'''}), as opposed to solutions
supported along a line.
\end{remark}

\begin{remark}
The ideas above can be used to provide an analogous combinatorial
description for the supports of certain solutions of $H_A(I,\beta)$
when $I$ is a general $\ZZ$-graded binomial ideal.  The key
observation is that if $\cC_{\rho,J}$ is a toral primary ideal, and
the parameter $\beta \in \CC^d$ is very generic inside
$\qdeg(\CC[\ttt]/\cC_{\rho,J})$, then the solutions of
$H_A(\cC_{\rho,J},\beta)$ are supported on translates of the
$L$-bounded components, where $L \subseteq \ZZ^J$ is the underlying
lattice of~$\rho$.  When $\cC_{\rho,J}$ is a primary component of~$I$,
this allows us to assert a lower bound on the number of series
solutions of $H_A(I,\beta)$ with the desired support.  Care must be
taken because $\qdeg(\CC[\ttt]/\cC_{\rho,J})$ could be partially or
entirely contained in the jump arrangement (Definition~\ref{d:jump}),
or $I_{\rho,J}$ could be an embedded prime, in which case the rank
at~$\beta$ need not equal a sum of multiplicities times volumes.
\end{remark}


\raggedbottom

\begin{thebibliography}{MMW05}

\bibitem[Ado94]{Adolphson}
Alan Adolphson, \emph{Hypergeometric functions and rings generated by
  monomials}, Duke Math. J. \textbf{73} (1994), no.~2, 269--290.

\bibitem[Ado99]{Adolphson-Rend}
Alan Adolphson, \emph{Higher solutions of hypergeometric systems and
  {D}work cohomology}, Rend. Sem. Mat. Univ. Padova \textbf{101}
  (1999), 179--190.  

\bibitem[App1880]{appell}
Paul Appell, \emph{Sur les s\'eries hyperg\'eometriques de deux
  variables et sur des \'equations diff\'erentielles lin\'eaires aux
  d\'eriv\'ees partielles}, Comptes Rendus \textbf{90} (1880),
  296--298.

\bibitem[BvS95]{bvs}
Victor V. Batyrev and Duco van Straten, 
\emph{Generalized Hypergeometric Functions and Rational Curves on
  Calabi-Yau Complete Intersections in Toric Varieties},
  Commun. Math. Phys. \textbf{168} (1995). 493--533.

\bibitem[Bor87]{borel}
Armand Borel, \emph{Algebraic $D$-modules}, Perspectives in
  Mathematics, 2. Academic Press, Inc., Boston, MA, 1987.

\bibitem[Bj{\"o}93]{Bjork2}
Jan-Erik Bj{\"o}rk, \emph{Analytic ${\mathcal D}$-modules and
  applications}, Mathematics and its Applications, vol.~247, Kluwer
  Academic Publishers Group, Dordrecht, 1993.

\bibitem[Bj{\"o}79]{Bjork}
Jan Erik Bj{\"o}rk, \emph{Rings of differential operators},
  North-Holland Mathematical Library, vol.~21, North-Holland
  Publishing Co., Amsterdam, 1979.

\bibitem[BH93]{BH93}
Winfried Bruns and J{\"u}rgen Herzog, \emph{Cohen-{M}acaulay rings},
  Cambridge Studies in Advanced Mathematics, vol.~39, Cambridge
  University Press, Cambridge, 1993.

\bibitem[CD07]{cattani-dick:binomial-ci}
Eduardo Cattani and Alicia Dickenstein, \emph{Counting solutions to
  binomial complete intersections}, J. of Complexity, Vol. \textbf{23}, 
Issue 1, Feb. 2007, 82--107.

\bibitem[CDL77]{quantum-mechanics}
Claude Cohen-Tannoudji, Bernard Diu and Franck Lalo\"{e},
  \emph{Quantum Mechanics, Volume One}, Wiley-Interscience, John Wiley
  \& Sons, New York, NY, 1977.

\bibitem[Cou95]{coutinho}
S. C. Coutinho, \emph{A primer of algebraic $D$-modules}, London
  Mathematical Society Student Texts, vol. 33, Cambridge University
  Press, Cambridge, 1995.

\bibitem[DMM08]{primDecomp}
Alicia Dickenstein, Laura~Felicia Matusevich, and  Ezra Miller,
  \emph{Combinatorics of binomial primary decomposition}, 
{math.AC/08033846}.

\bibitem[DMS05]{dms}
Alicia Dickenstein, Laura~Felicia Matusevich, and Timur Sadykov,
  \emph{Bivariate hypergeometric {$D$}-modules},
  Adv. Math. \textbf{196} (2005), no.~1, 78--123.

\bibitem[DS07]{ds06}
Alicia Dickenstein and Timur Sadykov, \emph{Bases in the solution
  space of the Mellin system}, Mat. Sbornik \textbf{198:9} (2007), 59-80.

\bibitem[Erd50]{erdelyi}
Arthur Erd\'elyi, \emph{Hypergeometric functions of two variables},
  Acta Math. \textbf{83} (1950), 131--164.

\bibitem[ES96]{binomialideals}
David Eisenbud and Bernd Sturmfels, \emph{Binomial ideals}, Duke
   Math. J. \textbf{84} (1996), no.~1, 1--45.

\bibitem[FS96]{fischer-shapiro}
Klaus~G. Fischer and Jay Shapiro, \emph{Mixed matrices and binomial
  ideals}, J. Pure Appl. Algebra \textbf{113} (1996), no.~1, 39--54.
 
\bibitem[Ful93]{fulton-toric}
William Fulton, \emph{Introduction to toric varieties}, Annals of
  Mathematics Studies, vol. 131, Princeton University Press,
  Princeton, NJ, 1993.
  
\bibitem[GGR92]{GGR}
I.~M. Gel{\cprime}fand, M.~I. Graev, and V.~S. Retakh, \emph{General
  hypergeometric systems of equations and series of hypergeometric
  type}, Uspekhi Mat. Nauk \textbf{47} (1992), no.~4(286), 3--82, 235.

\bibitem[GGZ87]{GGZ}
I.~M. Gel{\cprime}fand, M.~I. Graev, and A.~V. Zelevinski\u{\i},
  \emph{Holonomic systems of equations and series of hypergeometric
  type}, Dokl. Akad. Nauk SSSR \textbf{295} (1987), no.~1, 14--19.

\bibitem[GKZ89]{GKZ}
I.~M. Gel{\cprime}fand, A.~V. Zelevinski\u{\i}, and M.~M. Kapranov,
  \emph{Hypergeometric functions and toric varieties},
  Funktsional. Anal. i Prilozhen. \textbf{23} (1989), no.~2, 12--26.
  Correction in ibid, \textbf{27} (1993), no.~4,~91. 

\bibitem[GKZ94]{gkz-book}
I.~M. Gel{\cprime}fand, M.~M. Kapranov, and A.~V. Zelevinsky,
  \emph{Discriminants, resultants and multidimensional determinants},
  Mathematics: Theory \& Applications, Brikh\"auser Boston Inc.,
  Boston, MA, 1994.

\bibitem[God81]{godsil}
C.~D. Godsil, \emph{Hermite polynomials and a duality relation for matching
  polynomials}, Combinatorica \textbf{1} (1981), no.~3, 257--262.

\bibitem[GM92]{gm}
John P.~C. Greenlees and J.~Peter May, \emph{Derived functors of
  ${I}$-adic completion and local homology}, J. Algebra \textbf{149}
  (1992), no.~2, 438--453.

\bibitem[Ho99]{horja}
R.~P.~Horja, \emph{Hypergeometric functions and mirror symmetry in
  toric varieties}, math.AG/9912109.

\bibitem[Hor1889]{horn89}
J.~Horn, \emph{ {\"Uber die konvergenz der hypergeometrischen {R}eihen
  zweier und dreier Ver\"anderlichen}}, Math. Ann. \textbf{34} (1889),
  544--600.

\bibitem[Hor31]{horn31}
J.~Horn, \emph{Hypergeometrische {F}unktionen zweier
  Ver\"anderlichen}, Math. Ann. \textbf{105} (1931), no.~1, 381--407.

\bibitem[Hos06]{hosono}
Shinobu Hosono, \emph{Central charges, symplectic forms, and
  hypergeometric series in local mirror symmetry}, in: Mirror Symmetry V,
  AMS/IP Studies in Advanced Mathematics, Vol. 38, 2006, 405--440.

\bibitem[HLY96]{hly}
S.~Hosono, B.~H.~Lian, and S.~T.~Yau, \emph{GKZ-Generalized
  Hypergeometric Systems in Mirror Symmetry of Calabi-Yau
  Hypersurfaces}, Comm. Math. Phys. \textbf{182} (1996), no. 3,
  535--577.

\bibitem[Hot91]{equivariant}
Ryoshi Hotta, \emph{Equivariant {$D$}-modules}, math.RT/980502.

\bibitem[HS00]{hostenshapiro}
Serkan Ho{\c{s}}ten and Jay Shapiro, \emph{Primary decomposition of
  lattice basis ideals}, J. Symbolic Comput. \textbf{29} (2000),
  no.~4-5, 625--639, Symbolic computation in algebra, analysis, and
  geometry (Berkeley, CA, 1998).

\bibitem[It{\^o}51]{ito}
Kiyosi It{\^o}, \emph{Multiple {W}iener integral}, J. Math. Soc. Japan
  \textbf{3} (1951), 157--169.  

\bibitem[Kum1836]{kummer}
Ernst~Eduard Kummer, \emph{{\"U}ber die hypergeometrische {R}eihe
  {${F}(\alpha,\beta,x)$}}, J. f\"ur Math. \textbf{15} (1836).

\bibitem[Mel21]{mellin}
Hjalmar Mellin, \emph{R\'esolution de l'\'equation alg\'ebrique
  g\'en\'erale \`a l'aide de la fonction {$\Gamma$}},
  C.R. Acad. Sc. \textbf{172} (1921), 658--661.

\bibitem[Mil02b]{gradedGM}
Ezra Miller, \emph{Graded Greenlees--May duality and the \v Cech
  hull}, Local cohomology and its applications (Guanajuato, 1999),
  Lecture Notes in Pure and Appl. Math., vol. 226, Dekker, New York,
  2002, pp.~\mbox{233--253}.

\bibitem[MMW05]{MMW}
Laura~Felicia Matusevich, Ezra Miller, and Uli Walther,
  \emph{Homological methods for hypergeometric families},
  J. Amer. Math. Soc. \textbf{18} (2005), no.~4, 919--941.

\bibitem[MS05]{cca}
Ezra Miller and Bernd Sturmfels, \emph{Combinatorial commutative
  algebra}, Graduate Texts in Mathematics, vol. 227, Springer-Verlag,
  New York, 2005.

\bibitem[Oku06]{okuyama}
Go Okuyama, \emph{A-Hypergeometric ranks for toric threefolds},
  Internat. Math. Res. Notices \textbf{2006}, Article ID 70814,
  38~pages.

\bibitem[PWZ96]{A=B}
Marko Petkov{\v{s}}ek, Herbert~S. Wilf, and Doron Zeilberger,
  \emph{${A}={B}$}, A K Peters Ltd., Wellesley, MA, 1996.

\bibitem[Sad02]{timur}
Timur Sadykov, \emph{On the Horn system of partial differential
  equations and series of hypergeometric type},
  Math. Scand. \textbf{91} (2002), no.~1, 127--149.

\bibitem[SB02]{numerical-analysis-book}
J.~Stoer and R.~Bulirsch, \emph{Introduction to numerical analysis},
  third ed., Texts in Applied Mathematics, vol.~12, Springer-Verlag,
  New York, 2002.

\bibitem[SK85]{multiple-gauss-book}
H.~M. Srivastava and Per~W. Karlsson, \emph{Multiple {G}aussian
  hypergeometric series}, Ellis Horwood Series: Mathematics and its
  Applications, Ellis Horwood Ltd., Chichester, 1985.

\bibitem[SST00]{SST}
Mutsumi Saito, Bernd Sturmfels, and Nobuki Takayama, \emph{Gr\"obner
  {D}eformations of {H}ypergeometric {D}ifferential {E}quations},
  Springer-Verlag, Berlin, 2000.

\bibitem[SW08]{uli06}
Mathias Schulze and Uli Walther, \emph{Irregularity of hypergeometric systems 
via slopes along coordinate subspaces}, to appear: Duke Math. J.

\bibitem[Sla66]{slater}
Lucy Joan Slater, \emph{Generalized hypergeometric functions},
  Cambridge University Press, 1966.

\bibitem[Stu00]{Sturmfels:alg-hyp}
Bernd Sturmfels, \emph{Solving algebraic equations in terms of
  {$\mathcal{A}$}-hypergeometric series}, Discrete Math. \textbf{210}
  (2000), no.~1-3, 171--181, Formal power series and algebraic
  combinatorics (Minneapolis, MN, 1996).


\end{thebibliography}
\def\cprime{$'$} \def\cprime{$'$}
\providecommand{\MR}{\relax\ifhmode\unskip\space\fi MR }
\providecommand{\MRhref}[2]{%
  \href{http://www.ams.org/mathscinet-getitem?mr=#1}{#2}
}
\providecommand{\href}[2]{#2}

\end{document}